\documentclass[12pt,a4paper]{article}
\usepackage[latin1]{inputenc}
\usepackage{amsmath}
\usepackage{amsfonts}
\usepackage{amssymb}
\usepackage{graphicx}
\usepackage{color}
\usepackage{marginnote}
\usepackage{algorithmic}
\usepackage{algorithm}
\usepackage{geometry}
\def\cD{{\cal{D}}}
\def\S{{\bf S}}
\def\bZ{\mathbb{Z}}
\makeatletter      
\renewcommand{\ps@plain}{%
	\renewcommand{\@oddhead}{\textrm{A data driven Koopman-Schur decomposition }\today\hfil\textrm{\thepage}}%
	\renewcommand{\@evenhead}{\@oddhead}%
	\renewcommand{\@oddfoot}{}
	\renewcommand{\@evenfoot}{\@oddfoot}}
\makeatother     

 \newtheorem{theorem}{Theorem}[section]
 \newtheorem{proposition}{Proposition}[section]
 \newtheorem{lemma}{Lemma}[section]
\newtheorem{definition}{Definition}[section]
 \newtheorem{corollary}[theorem]{Corollary}
\newtheorem{example}[theorem]{Example}
\newtheorem{remark}[theorem]{Remark}

\newcommand{\DDS}{\mathbf{F}}

\newcommand{\eigA}{G}

\newcommand{\ecv}{\mathfrak{e}}

\newcommand{\X}{\mathbf{X}}
\newcommand{\Y}{\mathbf{Y}}

\newcommand{\x}{\mathbf{x}}
\newcommand{\y}{\mathbf{y}}
\newcommand{\w}{\mathbf{w}}
\newcommand{\kb}{\mathbf{k}}
\newcommand{\bfz}{\mathbf{z}}
\newcommand{\bfu}{\mathbf{u}}

\newcommand{\PR}{\mathbf{P}}

\newcommand{\ones}{\vec{\mathbf{1}}}

\newcommand{\bfalpha}{{\boldsymbol\alpha}}

\newcommand{\bfzeta}{{\boldsymbol\zeta}}
\newcommand{\bfdelta}{{\boldsymbol\delta}}
\newcommand{\bfrho}{{\boldsymbol\rho}}

\newcommand{\bfpsi}{{\boldsymbol\psi}}

\newcommand{\bfPhi}{{\boldsymbol\Phi}}
\newcommand{\bfPsi}{{\boldsymbol\Psi}}

\newcommand{\bfun}{{\boldsymbol \psi}} 
\newcommand{\bfunv}{{\boldsymbol\psi}}
\newcommand{\kerf}{{\boldsymbol f}}
\newcommand{\fspace}{{\mathcal{F}}}
\newcommand{\eigf}{{\boldsymbol\phi}}
\newcommand{\Pond}{{\boldsymbol \Omega}}
\newcommand{\pond}{{\boldsymbol \omega}}

\newcommand\restr[2]{{
		\left.\kern-\nulldelimiterspace 
		#1 
		\vphantom{\big|} 
		\right|_{#2} 
	}}
		
\newcommand{\R}{ \ensuremath{\mathbb{R}}}  
\newcommand{\Cz}{ \ensuremath{\mathbb{C}}}  
\newcommand{\0}{\mathbf{0}}
\newcommand{\Id}{\mathbb I}
\newcommand{\KO}{\mathbb U}

\def\cD{{\cal{D}}}
\newcommand{\cH}{\mathcal{H}}

\newcommand{\cB}{\mathcal{B}}
\newcommand{\caA}{\mathcal{A}}
\def\S{{\bf S}}
\def\bZ{\mathbb{Z}}
\def\cB{{\cal{B}}}

\begin{document}
\title{A data driven Koopman-Schur decomposition for computational analysis of nonlinear dynamics\thanks{
This material is based upon work supported by the DARPA Small Business Innovation Research Program (SBIR) Program Office under Contract No. W31P4Q-21-C-0007 to AIMdyn, Inc. Any opinions, findings and conclusions or recommendations expressed in this material are those of the authors and do not necessarily reflect the views of the DARPA SBIR Program Office. \textbf{Distribution Statement A:  Approved for Public Release, Distribution Unlimited.} Submitted to SIAM Journal on Scientific Computing.
}}
\author{Zlatko Drma\v{c}\thanks{University of Zagreb, Faculty of Science, Department of Mathematics, Bijeni\v{c}ka cesta 30, 10000 Zagreb, Croatia. The author is also supported by the Croatian Science Foundation (CSF) grant IP-2019-04-6268 \emph{``Randomized low rank algorithms and applications to parameter dependent problems''}.} \and Igor Mezi\'{c}\thanks{University of California at Santa Barbara and AIMdyn, Inc., 1919 State St., Ste. 207, Santa Barbara, CA 93101. The author is also supported by the AFOSR Award FA9550-22-1-0531 and the ONR Award N000142112384. }}

\date{}

\maketitle

\begin{abstract}
This paper introduces a new theoretical and computational framework for a data driven Koopman mode analysis of nonlinear dynamics. To alleviate the potential problem of ill-conditioned eigenvectors in the existing implementations of the Dynamic Mode Decomposition (DMD) and the Extended Dynamic Mode Decomposition  (EDMD), the new method introduces a Koopman-Schur decomposition that is entirely based on unitary transformations. The analysis in terms of the eigenvectors as modes of a Koopman operator compression  is replaced with a modal decomposition in terms of a flag of invariant subspaces that correspond to selected eigenvalues. The main computational tool from the numerical linear algebra is the partial ordered Schur decomposition that provides convenient orthonormal bases for these subspaces. In the case of real data, a real Schur form  is used and the computation is based on real orthogonal transformations. 
The new computational scheme is presented in the framework of the Extended DMD and the kernel trick is used. %
\end{abstract}

\section{Introduction and preliminaries}\label{S=Intro}
High resolution measurements (using e.g. particle image velocimetry (PIV), high speed cameras, magnetometers, microwave radiometers) or high fidelity numerical simulations of the dynamics under study (given e.g. as a system of ordinary differential equations obtained in a semi-discretization of partial differential equations) often provide an abundance of data.  After preprocessing and formatting, the data is represented as a discrete time sequence of data snapshots $\x_1, \x_2, \ldots$ in  the state space that is set to be $\Cz^n$ or $\R^n$, and assumed to be driven by a mapping $\x_{i+1}=\DDS(\x_i)$, for some initial state $\x_1$ and with given time step $\delta t$. The mapping $\DDS$ may not be known; in a numerical simulation it represents software implementation of a numerical algorithm. The dimension $n$ ranges, depending on the application, from moderate to extremely large. For instance, if the sequence $\x_1, \x_2, \ldots$ is obtained by discretization with time step $\delta t$ of a time dependent partial differential equation, then $\x_i$ is represents the solution on a spatially discretized domain at the discrete time moment $i$ and the dimension $n$ (the number of grid points in the discretized domain) can be in hundreds of thousands or in millions. In such cases the natural structure of the data is a fourth order tensor (three space dimensions times the time), and the ambient space $\Cz^n$ or $\R^n$ is given implicitly by a reshaping isomorphism (vectorization).

\vspace{-3mm}

\paragraph{The Dynamic Mode Decomposition (DMD).} 
Revealing relevant features, latent in the data but valuable for understanding and control of the underlying processes, requires suitable data mining techniques. 
The Dynamic Mode Decomposition (DMD) \cite{Schmid-DMD-2208} is such a tool of the trade. It has been initially motivated by the problems in the computational fluid dynamics, but its success ranges across a wide spectrum of applications, such as combustion, \cite{Combustion-Inst-Flame-Images-2016}, robotics \cite{dmd-robotics-2015}, power grid \cite{susuki2016applied} and automatic control \cite{mauroy2020koopman}. 

The tacit assumption in the DMD is that the given sequence of  snapshots $\x_1, \ldots, \x_m, \x_{m+1}$ is generated by a linear operator (matrix) $A_n$ such that, given $\x_1$, $\x_{i+1}=A_n \x_i$ for $i=1,\ldots,m$. If we set $\X=(\x_1,\ldots,\x_m)$, $\Y=(\x_2,\ldots,\x_{m+1})$, then 
$\Y = A_n\X$. Since the action of $A_n$ on $\X$ is given, we can use the Rayleigh-Ritz procedure to extract from the range of $\X$  approximate eigenpairs of $A_n$. To that end, a POD basis is computed using the SVD of $\X$, $\X = U_r \Sigma_r V_r^*$, where $r$ is the rank of $\X$. In practice, this will be the numerical rank, i.e. small singular values of $\X$ will be truncated. Using the orthonormal basis $U_r$, the Rayleigh quotient is computed as $A=U_r^* A_n U_r = U_r^* \Y V_r\Sigma_r^{-1}$. If we diagonalize $A$, $A=\eigA \Lambda \eigA^{-1}$ with $\Lambda=\mathrm{diag}(\lambda_i)_{i=1}^r$, then the Ritz pairs of $A_n$ are $(U_r \eigA, \Lambda)$. If we set $K=U_r \eigA=(\kb_1,\ldots,\kb_r)$, then $A_n \kb_i \approx \lambda_i\kb_i$. Note that there is no need to form the matrix $A_n$ explicitly. A software implementation of the DMD has been recently included in the LAPACK library \cite{LAPACK}, \cite{lawn-298-drmac}, \cite{lawn-300-drmac}.

Now, revealing the structure of the snapshots (e.g. discovering coherent structures in the fluid flow) is based on a spatio-temporal representation of the snapshots using selected $\ell\leq r$ modes $\kb_{j_1}, \ldots , \kb_{j_\ell}$: For given weights $\omega_i\geq 0$, the task is to determine coefficients $\alpha_1,\ldots,\alpha_\ell$ so that
\begin{equation}\label{eq:KMD-1}
\x_i \approx \sum_{k=1}^\ell \kb_{j_k} \alpha_k \lambda_{j_k}^{i-1}\;\;\mbox{and}\;\; \sum_{i=1}^{m+1}\omega_i^2 \|\x_i - \sum_{k=1}^\ell \kb_{j_k} \alpha_k \lambda_{j_k}^{i-1}\|^2 \longrightarrow \min_{\alpha_k} .
\end{equation}
Note that in the above formula, the representation of the snapshot $\x_{i+1}=A_n\x_i$ is obtained simply by raising by one the powers $\lambda_{j_k}^{i-1}$ in the representation of $\x_i$. The key is that the coefficients in the representation  (\ref{eq:KMD-1}) of $\x_i$ also contain the time stamp. This is a basis for forecasting -- extrapolation into the future by raising the powers of the $\lambda_{j_k}$'s in the representation (\ref{eq:KMD-1}) of the last available (i.e. the present) snapshot. The representation of the snapshots (\ref{eq:KMD-1}) in terms of selected modes $\kb_{j}$ is most dynamics revealing if it succeeds with small $\ell$, and such a sparse representation can be obtained using a sparsity promoting DMD   \cite{jovschnicPOF14}.

For more details on the variations of the DMD and its various applications see e.g.
\cite{schmid2011}, \cite{Chen:2012jh}, 
\cite{2015arXiv150203854H}, \cite{Dawson2016}, 
\cite{Hemati:2014jm}, \cite{Takeishi-PhysRevE.96.033310}, \cite{Takeishi-ijcai2017-392}, \cite{Takeishi-8296769}, \cite{Takeishi-NIPS2017_6713}, \cite{dmd-control-doi:10.1137/15M1013857}, \cite{varpro-opt-dmd-doi:10.1137/M1124176}, \cite{book-Koopman-Sys-Control}.

\paragraph{The Koopman connection and the Extended DMD (EDMD).} The theoretical underpinning of the DMD is in  close connection to the Koopman (composition) operator $\KO$ associated with the discrete dynamical system $\x_{i+1}=\DDS(\x_i)$. In fact, the leftmost formula in (\ref{eq:KMD-1}) is a version of the general Koopman Mode Decomposition, first derived in \cite{mezic2005spectral}. $\KO$ is defined on a class of functions $\mathcal{F}$ equipped with a Hilbert space structure (typically $L^2$ for measure-preserving systems, see \cite{mezic2020spectrum} for the Hilbert space construction in the dissipative case) as 
$\KO f = f\circ\DDS$, $f\in\mathcal{F}$. {Here $\circ$ denotes the composition of mappings. The functions $f\in\mathcal{F}$ are called the observables of the system.}  Note that $(\KO f)(\x_i)=f(\DDS(\x_i))=f(\x_{i+1})$, so that $\KO$ acts similarly as the matrix $A_n$ in the DMD, but in a more general setting. The introduction of $\KO$ is actually a linearization of $\DDS$, albeit an infinitely dimensional one. 
The connection with the DMD  is explicitly stated in the Extended DMD \cite{williams-2015-EDMD}, which
encompasses the DMD as a special case, and can be understood as a Galerkin type approximation of the eigenvalues and eigenfunctions of the Koopman operator, based on its compression to an $N$ dimensional subspace of $\mathcal{F}$. 
For an excellent review see e.g. \cite{Brunton-etal-SIREW-2022}.

It is known that, under certain assumptions, when the number of snapshots in the limit  $m\rightarrow\infty$ and the dimension $N\rightarrow\infty$, the EDMD in several  ways (including the convergence of the eigenfunctions and the eigenvalues along subsequences) converges to the operator $\KO$, see  \cite{korda-mezic-EDMD-convergence}. The use of Hankel-DMD \cite{arbabi2017ergodic} under the assumption of existence of finite-dimensional invariant subspace leads to stronger results, that were extended to convergence in pseudospectral sense in \cite{mezic2022numerical}, without the assumption on existence of a finite-dimensional invariant subspaces, under the assumption of discrete spectrum of the Koopman operator.

Here too, the eigenfunctions of $\KO$ are numerically approximated using a discrete compression of $\KO$ to an $N$-dimensional subspace $\mathcal{F}_N=\mathrm{span}(\bfun_1,\ldots,\bfun_N) \subset\mathcal{F}$, where $N\geq n$ and potentially $N\gg n$. Similarly as in the DMD, a Rayleigh quotient matrix $U_N$ is diagonalized and its eigenvectors are used to construct discrete approximations of certain eigenfunctions of $\KO$. In a special case, the matrix $U_N$ and the DMD matrix $A_n$ are transposes of each other,\footnote{See Section \ref{S=EDMD}.} and the eigenvectors of $A_n$ are the modes in the data snapshot representation (\ref{eq:KMD-1}).

\paragraph{A quandary about the eigenvectors.} There is one important detail that is, to the best of our knowledge, not properly addressed in the literature. If in an application we compute a spectral decomposition $A = \eigA\Lambda \eigA^{-1}$ of a highly non-normal matrix $A$, assuming $A$ is diagonalizable, then the matrix $\eigA$ of the eigenvectors can be severely ill-conditioned and in that case the matrix of the modes $K=U_r \eigA$ must also be ill-conditioned. Ill-conditioned means close to singularity and numerical computation (such as representing the snapshots as linear combinations of the eigenvectors (\ref{eq:KMD-1})) becomes prone to large errors.

Even if $A$ is not diagonalizable (i.e. it is defective), a numerical software will compute spectral information that corresponds to a nearby $A+\delta A$ (this is established by a backward or a mixed error analysis of the algorithm that is implemented in the software) that is most likely diagonalizable. (Recall  that non-diagonalizable matrices are a nowhere dense subset - the diagonalizable ones are an open dense set.)
As a result, a numerical software is likely to return a full set of eigenvectors, but alas, not of $A$. A simple example illustrates this. Let 
$
A = \left(\begin{smallmatrix} 1 & 1 \cr 0 & 1\end{smallmatrix}\right)$, $A + \delta A = \left(\begin{smallmatrix} 1 & 1 \cr \epsilon & 1\end{smallmatrix}\right)$, where $\epsilon >0$ is a small parameter. $A$ is a Jordan block with double eigenvalue $\lambda=1$ of geometric multiplicity one.
The eigenvalues of $A+\delta A$ are $\lambda_{1,2}= 1\pm\sqrt{\epsilon}$ and the spectral decomposition of $\widetilde{A}=A+\delta A$ reads
$$
\widetilde{A} = \widetilde{\eigA} \widetilde{\Lambda} \widetilde{\eigA}^{-1},\; \widetilde{\eigA} = \begin{pmatrix} 1 & 1 \cr \sqrt{\epsilon} & -\sqrt{\epsilon}\end{pmatrix}\! ,\;
\widetilde{\Lambda} = \begin{pmatrix} 1+\sqrt{\epsilon} & 0 \cr 0 & 1-\sqrt{\epsilon}\end{pmatrix}\! ,\;
\widetilde{\eigA}^{-1}=\frac{1}{2}\begin{pmatrix} 1 & 1/\sqrt{\epsilon}\cr 1 & -1/\sqrt{\epsilon}\end{pmatrix}.
$$
Note that $\|\widetilde{\eigA}\|_2\|\widetilde{\eigA}^{-1}\|_2 > 1/\sqrt{\epsilon}$ and that no scaling of the eigenvectors can improve this ill-conditioning with a factor better than $\sqrt{2}$. This is a consequence of a van der Sluis theorem \cite{slu-69} on optimally scaled matrices and the fact that both columns of $\widetilde{\eigA}$ are of Euclidean norm $\sqrt{1+\epsilon}\approx 1$. 

Hence,  theoretical considerations that are based on the nontrivial (i.e. non-diagonal) Jordan normal form and the generalized eigenvectors are not robust for practical computation -- an arbitrarily small perturbation can  change the Jordan structure entirely and e.g. make the previously identified generalized eigenvectors to simply disappear and replace them with potentially ill-conditioned eigenvectors.  
For a more general result on optimal scaling of the eigenvectors in the case of multiple eigenvalues and on the Jordan normal form  from a numerical analyst's perspective we recommend  \cite{demmel-thesis-83}.

{
The analysis of the case when the Koopman operator itself has a degeneracy that leads to generalized eigenfunctions was presented in \cite{mezic2020spectrum}.\footnote{ The notion of generalized eigenfunctions used in \cite{mezic2020spectrum} is in the same sense as generalized eigenvectors of a matrix - i.e.functions that span an  invariant subspace of the Koopman operator of dimension bigger than $1$ that contains only $1$ eigenfunction. This is in contrast with the notion of generalized eigenfunctions in functional analysis \cite{gel2016generalized} that are objects in Koopman theory associated with continuous spectrum. } A simple example in which generalized eigenfunctions arise is the linear system
\begin{equation}
    \dot \y=A\y,
\end{equation}
where $\y=(y_1,y_2)^T$, $A = \left(\begin{smallmatrix} 1 & 1 \cr 0 & 1\end{smallmatrix}\right)$. Since linear (generalized) eigenfunctions of the Koopman operator in the linear case are given by $\left<\y,\w \right>,$
where  $\w$ is a left (generalized) eigenvector of $A$ and $\langle\cdot,\cdot\rangle$ is the complex inner product defined by
$
\langle\x,\y\rangle=\sum_i \x_i \cdot \overline{\y_i}.
$
The eigenvector $\w_1$ associated with eigenvalue $1$ is $\w_1=(0,1)^T$, and the generalized eigenvector can be chosen to be $(1,0)^T$. Therefore $y_2$ is an eigenfunction, and $y_1$ a generalized eigenfunction of the system. 

The paper \cite{schmid2007nonmodal} contains examples of non-normal operators in fluid  flows, that has substantial impact on treatment of stability of such flows. The most well-known example of such non-normality is the formulation of the viscous stability problem for parallel shear flows, in which linearization of the Navier-Stokes equation leads to the Orr-Sommerfeld linear partial differential equation whose solutions exhibit non-normal behavior. Numerical simulations reveal that the computed eigenvectors are indeed highly ill--conditioned.}

\paragraph{An alternative proposed in this paper.} We explore a different method to numerically extract spectral information about $\KO$ and to deploy it in applications with the same functionality as provided by the DMD and the EDMD.  
The key tool in the new approach is the Schur triangular form, which is in the framework of the numerical linear algebra a viable alternative to the diagonalization $A=\eigA\Lambda \eigA^{-1}$ or computation of the Jordan canonical form.
\begin{theorem}\label{TM:Schur-form}
(Schur form)
Let $A\in\Cz^{r\times r}$ and let $\lambda_1,\ldots, \lambda_r$ be its eigenvalues, listed in an arbitrary order. Then there exist a unitary $Q$ and an upper triangular $T$ such that $T_{ii}=\lambda_i$, $i=1,\ldots, r$, and 
$
A = Q T Q^* .
$
\end{theorem}\label{remark-Schur-form}
\begin{remark}
{\em
If we write $T=\mathrm{diag}(\lambda_i) + N$, then $A = Q \mathrm{diag}(\lambda_i) Q^* + QNQ^*$ is the sum of a normal and a nilpotent matrix; $(QNQ^*)^r=\0$.
If we let $T_{1:i,1:i}$ denote the leading $i\times i$ submatrix of $T$, and $Q_{:,1:i}$ the submatrix of the leading $i$ columns of $Q$, then
$AQ_{:,1:i}=Q_{:,1:i}T_{1:i,1:i}$, i.e. the matrices $Q_{:,1:i}$ span a flag of invariant subspaces of $A$, with the corresponding compressions (Rayleigh quotients) given by $T_{1:i,1:i}=Q_{:,1:i}^* A Q_{:,i}$. The corresponding orthogonal projectors $P_i = Q_{:,1:i}Q_{:,1:i}^*$ are increasing in the Loewner partial order: $P_{i+1} \succeq P_i$, meaning that $P_{i+1}-P_i=Q_{:,i}Q_{:,i}^*$ is positive semidefinite.
}
\hfill $\boxtimes$\end{remark}
\begin{remark}
{\em
If $A=QTQ^*$, then $A^* = Q T^* Q^* = (Q\Pi) (\Pi^T T^* \Pi) (Q\Pi)^*$ is the Schur form of $A^*$, where $\Pi$ is the permutation matrix that reverses the order of columns, so that $\Pi^T T^* \Pi$ is upper triangular. Similarly, the Schur form of the transposed matrix reads
$A^T= (Q\Pi)^{*T} (\Pi^T T^T \Pi) (Q\Pi)^T$.
}
\hfill $\boxtimes$\end{remark}
An advantage of this decomposition is that it is numerically computed using only unitary transformations, and that it does not refer to the Jordan normal form of $A$.
Having $A = QTQ^*$, we have in the DMD $A_n (U_r Q) \approx (U_r Q) T$, and in the EDMD this yields a partial Schur decomposition of the matrix compression of $\KO$. 

When using the Schur form, in the general non-normal case we do not have direct access to the individual eigenvectors\footnote{The exception is the first column of $Q$, $A Q(:,1)=T_{11} Q(:,1)$.} (discretized eigenfunctions) anymore -- instead, we have a flag of nearly invariant subspaces and additional work is needed to get the same functionality as in the EDMD. Our proposed solution is based on the numerical linear algebra techniques, in particular on the ordered Schur form. 
A special case that covers Hermitian and skew-Hermitian DMD (thus diagonal Schur form) is recently published in \cite{lawn-300-drmac}. We have already used an ordered partial Schur decomposition\footnote{We called it partial Schur Modal Decomposition  (PSMD)} in \cite[Section 5.3]{DDMD-SISC-2018} to improve the Optimal Mode Decomposition \cite{wynn_pearson_ganapathisubramani_goulart_2013}, and we discussed some details of using a partial Schur decomposition in a DMD computational framework.
\footnote{
A similar idea is also used recently in \cite{9107086}, but only partially. The authors of \cite{9107086} seemed unaware of our 2018 paper \cite{DDMD-SISC-2018}.} 

For the behaviour of the computed matrix Schur forms in the limit for $m\rightarrow\infty$ and $N\rightarrow\infty$, we must consider a Schur form of $\KO$.  Since the matrix Schur form can be extended to spectral operators \cite{rogers1990triangular}, we proceed beyond the numerical procedure and complete the new proposed framework with the Schur decomposition of the Koopman operator and the convergence of the numerical method.

\paragraph{Outline of the paper.} The rest of the paper is organized as follows. In Section \ref{S=Koopman=Schur-operator} we set the operator theoretic stage for the Koopman operator and its Schur form -- from the definition of the Koopman semigroup, nonlinear dynamics, to properly defined function spaces.  

In Section \ref{S=EDMD}, we review the EDMD in detail. The Koopman operator based interpretation is stressed (as opposed to the mere regression) and the details of a data driven compression of the operator and its Rayleigh quotient are given. Further, the kernel trick is explained and elements of convergence in appropriate Hilbert spaces
are given. Together with details such as evaluation of the approximate eigenfunctions on the snapshots, this section serves to define the framework and set the stage for our new proposed Koopman-Schur  Subspace Modal Decomposition (KS-SSMD), that is described in detail in Section \ref{S=KSD}.     
We show that the KS-SSMD has all the functionalities of the EDMD, and, in particular,  that it can be used to reveal coherent structures (Section \ref{SS=Observ_Representation_Schur_f}) and for forecasting (Section \ref{SS=KS_forecasting}). 
The flexibility of working with only a subset of selected eigenvalues is achieved using the ordered Schur form, that is outlined in Section \ref{SS=ordered-Schur}. The new approach is put to test in Section \ref{S=NumericalEvaluation}. We use numerical examples to illustrate subtle numerical issues and the potential of the KS-SSMD, in particular good forecasting skills. 
This report concludes with Section \ref{S=Conclusion}, that contains a review of our current work on the further development of the method.

{
\section{Koopman operator and its Schur form}\label{S=Koopman=Schur-operator}

Let $T:M\rightarrow M$ be a dynamical system on a set $M$.  Let the Koopman operator $\KO$ associated with $T$,  be  defined by 
\begin{equation*}
    \KO f=f\circ T, f\in \cal F
\end{equation*}
where $\cal F$ is a space of observables closed under the action of $\KO$ i.e. $\KO(\cal F)\subset \cal F.$

We recall the definition of a spectral operator \cite{dundford:1954}:
\begin{definition}
  A bounded operator $S$ on a space $\cal F$ is of scalar type  if it satisfies \
\begin{equation}
S=\int_\Cz  \lambda dP_{\lambda},
\end{equation}
where $dP$ is the resolution of the identity for $S$.
An operator $U$ is {\it spectral},  if it can be decomposed as
$$
U=S+Q_0,
$$
where $S$ is of scalar type, and $Q_0$ is quasi-nilpotent, i.e.
$$
\lim_{n \rightarrow \infty}(Q_0^n)^{1/n}=0.
$$
\end{definition}
The following theorem is due to Rogers \cite{rogers1990triangular}.
\begin{theorem}\label{ther:Rodg} Let U be a spectral operator. Then $U= N + Q_1$
where $N$ is normal, $Q_1$ is quasinilpotent, and $\sigma(N) = \sigma(U)$. Also $Q_1N$, $NQ_1$, and
$NQ_1 - Q_1N$ are quasinilpotent. Furthermore there is a countable nest of selfadjoint
projections ${\cal N} = \{ P_\alpha \}_{\alpha \in \Gamma}$  with the following properties:
\begin{enumerate}
\item Each $P_\alpha$, reduces $N$ and is invariant under $Q$,
\item The operator $N$ is in the $C^*$-algebra generated by ${\cal N}$,
\item The nest ${\cal N}$ separates $\sigma(U)$ in the sense that given any $\lambda,\mu \in \sigma(U)$,
there exists   a $P_\alpha$, such that $\sigma(U_{P_\alpha{\cal F}})$ contains exactly one of $\lambda,\mu$.
\end{enumerate}
\end{theorem}

We work in continuous time $t$, where $\S^t:M\rightarrow M$ is a group of transformations of $M$. The type system we consider here is the one that has a global (Milnor) attractor $\caA$ of zero Lebesgue measure,  where $\caA$ is equipped with a physical measure \cite{climenhagaetal:2017}. The operator $U$ can be obtained from the family of operators $U^t$ defined by  $U^tf=f\circ \S^t$ by sampling the evolution of the system at discrete times $n\Delta t,n\in \bZ$, with $T=S^{\Delta t}$.
Let $\cD$ be a compact forward invariant set of $\S^t$, and  $\mu$ the physical invariant measure on the (Milnor) attractor $\caA \subset \cD$, with $\mu(\caA)=1$. Let $\cH_\caA$ be the space of square-integrable observables on $\caA$. Note that $U^t$ restricted to $\cH_\caA$ is unitary. Thus the theorem above yields $U^t|_{\cH_\caA}=N^t$, i.e. $U^t$ is a normal operator when restricted to $\cH_\caA$. The functions in $\cH_\caA$ can typically be thought of as being defined on the whole basin of attraction of  $\caA$.  Let $\tilde\cH_\cB\subset C(\cD)$ be a Hilbert space of continuous functions $f:\cD\rightarrow \mathbb{C}$ orthogonal 
to $\cH_\caA$ with respect to $\mu$, that satisfy
\begin{equation}
\int_\cD f\phi d\mu=0, \ \forall \phi \in \cH_\caA.
\label{eq:ort}
\end{equation}
We assume that the tensor product space  $\tilde\cH=\cH_\caA\otimes \tilde\cH_\cB$ is invariant under $U^t$.  It is clear that the space $\cH_\caA$ is itself invariant under $U^t$. So is $\tilde\cH_\cB$ \cite{mezic:2020}:
\begin{lemma} $\tilde\cH_\cB$ is invariant under $U^t$.
\end{lemma}
Define $\cH_\cB=\tilde \cH_\cB\cup {\mathbf 1}$, where ${\mathbf 1}$ is the constant unit function on $\cD$.
Define the tensor product space \begin{equation}\cH=\cH_\caA\otimes \cH_\cB.\end{equation} Clearly,
$U^t=U^t|_{\cH_\caA}\otimes U^t|_{\cH_\cB}$ on $\cH$. \hfill$\boxtimes$

A suitable space $\cH_\cB$ in this context is the so-called Modulated Fock Space \cite{mezic:2020} where principal eigenfunctions of the Koopman operator are used in the definition of the space. Using the techniques in \cite{mezic:2020} the following can be shown:
\begin{theorem}
Let $U$ be a discretization of the Koopman family $U^t$ for a system with an $n$-dimensional Milnor attractor $\caA$. Let $\cH_\cB$ be the Modulated Fock Space. We assume that $U^t|_{\cH_\cB}$ is spectral. Then $U$ admits the Schur form
\begin{equation}
    U=N+Q_1
\end{equation}
where $N$ is normal and $Q_1$ is quasi-nilpotent.
\end{theorem}
{\sc Proof:}
We know that $U^t|_{\cH_{\caA}}$ is scalar, as it is the space of square integrable functions with respect to $\mu$ rendering $U^t|_{\cH_\caA}$ unitary. The fact that  $U^t|_{\cH_\cB}$ is spectral immediately renders  $U^t|_{\cH_\caA}\otimes U^t|_{\cH_\cB}$ spectral, and therefore Theorem \ref{ther:Rodg} applies.
\ \hfill$\boxtimes$
}
\section{Extended DMD}\label{S=EDMD}
In this section, we review the EDMD following the seminal work \cite{williams-2015-EDMD}, upon which we build our discrete Koopman-Schur decomposition. For the sake of completeness and for the readers' convenience when comparing the proposed method with the EDMD, we provide full description of the main ingredients of the EDMD. For further study of the EDMD the reader is referred to \cite{williams-2015-EDMD}, and for convergence analysis we refer to \cite{korda-mezic-EDMD-convergence}, \cite{Colbrook-Townsend-2021}, \cite{colbrook_ayton_szoke_2023}.

In Section \ref{SS=Num-Compress} we first describe the construction of a finite dimensional data driven  compression $U_N$ of the Koopman operator to a selected $N$-dimensional function subspace $\fspace_N$, and in Section \ref{SS=Convergence} we analyze the convergence when $N\rightarrow\infty$ and the number of available snapshot pairs $m\rightarrow\infty$. If $m<N$, then the proper
method to approximate the eigenvalues and eigenvectors is the Rayleigh-Ritz extraction and in Section \ref{SS=RQ_UN} we review the details.
The kernel trick, that allows for an efficient computation with extremely large dimension $N$ is briefly described in Section \ref{SS=EDMD-kernek-trick}. The technical details of point-wise evaluation of the approximate eigenfunctions and representation of the data snapshots as linear combinations of the computed eigenfunctions are provided in Sections \ref{SS=evaluate-fs} and \ref{SS=EDMD-snap-repr}, respectively.
In addition, as a preparation and a motivation for the new method in Section \ref{S=KSD}, we discuss some numerical/implementation issues that are not tackled in  \cite{williams-2015-EDMD}.

%
%
\subsection{Numerical data driven compression of $\KO$}\label{SS=Num-Compress}
Numerical approximation of some eigenvalues and eigenfunctions of $\KO$ begins with a selection of a finite dimensional subspace
$
\fspace_N = \mathrm{Span}(\bfun_1,\ldots, \bfun_N) , 
$
where $\bfun_1,\ldots, \bfun_N$ are suitable scalar basis functions defined on $\Cz^n$ (or $\R^n$ for real data). In an ideal situation, an abundance of snapshot data $(\x_1,\y_1),\ldots (\x_m, \y_m)$ is available ($\y_i=\DDS(\x_i)$ and $m$ is as large as we can afford) and with a suitable choice of basis functions we can let $N\rightarrow\infty$ (i.e. we can make $N$ as large as desired -- if we take polynomials in $n$ variables of the total degree of at most $p$, then $N=\binom{n+p}{n}$). In fact, using a kernel trick with appropriate functions, we can have even (implicitly) $N=\infty$. Furthermore, the $\bfun_i$'s are preferably chosen from a sufficiently rich class of functions with good approximation properties in the ambient function space. We defer the convergence analysis (when $m, N\rightarrow\infty$) to 
Section \ref{SS=Convergence}, and proceed with a review of numerical details of matrix representation of the data driven compression, using material from \cite{math9172075}, \cite{lawn-298-drmac}, \cite{lawn-300-drmac}, \cite{Mezic-etall-COVID-2023}.  

Data driven setting means that we are supplied with data snapshots pairs $(\x_1,\y_1=\DDS(\x_1))$, $\ldots,$ $(\x_m, \y_m=\DDS(\x_m))$, and the values of the basis functions evaluated at these snapshots are conveniently tabulated in the $m\times N$ matrices $\bfPsi_x$, $\bfPsi_y$ as 
\begin{equation}\label{eq:data_matrices}
\bfPsi_x = \begin{pmatrix} \bfun_1(\x_1)  & \ldots & \bfun_N(\x_1) \cr\vdots & \ldots & \vdots \cr
\bfun_1(\x_m)  & \ldots & \bfun_N(\x_m)\end{pmatrix},\;\;
\bfPsi_y = \begin{pmatrix} \bfun_1(\y_1)  & \ldots & \bfun_N(\y_1) \cr\vdots & \ldots & \vdots \cr
\bfun_1(\y_m)  & \ldots & \bfun_N(\y_m)\end{pmatrix} .
\end{equation}
If the data are collected from a single trajectory, then $\y_i=\x_{i+1}$. In general, the data may be collected from several trajectories and the basic assumption is that the data is organized so that $\y_i=\DDS(\x_i)$.
It is reasonable to assume that $\bfun_1, \ldots, \bfun_N$ are linearly independent. However, if those functions are visible only on the snapshots $\x_i$, and $m < N$, then at most $m$ of them will be seen as independent -- in other words, the ranks of the matrices $\bfPsi_x$, $\bfPsi_y$ are at most $m$. In a data rich scenario with $m > N$, it is reasonable to assume that the rank is $N$, but a numericall algorithm should check the numerical rank in case of ill-conditioning.

Once the data matrices (\ref{eq:data_matrices}) are assembled, an operator compression $\PR_{\fspace_N}\restr{\KO}{\fspace_N} : \mathcal{F}_N \longrightarrow \mathcal{F}_N$ is computed numerically, where $\PR_{\fspace_N}$ is an appropriately defined projection onto $\fspace_N$. Since $\fspace_N$ is unlikely to be $\KO$-invariant, we have 
$$
\KO \bfun_i = \bfun_i\circ \DDS=\sum_{j=1}^N \bfu_{ji}\bfun_j + \rho_i,
$$
where $\bfu_{ji}$, $j=1,\ldots, N$,  are the coefficients that minimize the residual function $\rho_i$, $\int |\rho_i|^2 d\mu\rightarrow\min$. 
Since in a data driven setting the only available information is (\ref{eq:data_matrices}), the only feasible action is to solve the discrete least squares (LS) problems 
\begin{equation}\label{eq:LS-Ui}
\int \left| \sum_{j=1}^N \bfu_{ji}\bfun_j - \bfun_i\circ\DDS\right|^2  d\bfdelta_{m} 
\longrightarrow \min_{\bfu_{1i},\ldots,\bfu_{Ni}},\;\; i =1, \ldots, N, 
\end{equation}
i.e. with respect  to the empirical measure defined using the sum of the Dirac measures concentrated at the $\x_k$'s, $\bfdelta_{m} = (1/m)\sum_{k=1}^m\bfdelta_{\x_k}$. 
\begin{remark}
{\em 
In other words, (\ref{eq:LS-Ui}) computes the orthogononal projection of $\KO\bfun_i$ onto $\mathcal{F}_N$, computed in the space $L^2(\bfdelta_m)$. In the continuous setting in $L^2(\mu)$, it simply reads
$\int \left| \sum_{j=1}^N \bfu_{ji}\bfun_j - \bfun_i\circ\DDS\right|^2  d\mu 
\longrightarrow \min_{\bfu_{1i},\ldots,\bfu_{Ni}},\;\; i =1, \ldots, N.$ 
With properly weighted (quadrature based) algebraic least squares  rigorous error bounds are possible, see e.g. 
\cite{Colbrook-Townsend-2021}, \cite{colbrook_ayton_szoke_2023}.
}
\hfill $\boxtimes$\end{remark}

Since $(\KO\bfun_i)(\x_j)=(\bfun_i\circ \DDS)(\x_j)=\bfun_i(\y_j)$,  the minimization problem (\ref{eq:LS-Ui}) reads
\begin{equation}\label{eq:LS-Ui-2}
\frac{1}{m} \sum_{j=1}^m |\rho_i(\x_j)|^2 = 
\frac{1}{m}\left\|\! \left[
\bfPsi_x
\left(\begin{smallmatrix} \bfu_{1i}\cr \vdots \cr \bfu_{Ni}\end{smallmatrix}\right)\! - \! 
\left(\begin{smallmatrix} \bfun_i(\y_1)) \cr
\vdots  \cr
\bfun_i(\y_m) \end{smallmatrix}\right)\!
\right]\!\!\right\|_2^2 \longrightarrow \min_{\bfu_{1i},\ldots,\bfu_{Ni}},\;\; i=1,\ldots, N.
\end{equation}
To ensure well defined unique matrix representation of $\PR_{\fspace_N}{\KO}_{|\fspace_N}$ (independent of the rank of $\bfPsi_x$), an additional constraint $\sum_{j=1}^N |\bfu_{ji}|^2 \rightarrow\min$ is imposed for all $i$, yielding $(\bfu_{ji})_{j=1}^N = \bfPsi_x^\dagger (\bfun_i(\y_k))_{k=1}^m=\bfPsi_x^\dagger \bfPsi_y(;,i)$.
Here $\bfPsi_x^\dagger$ denotes the Moore-Penrose generalized inverse of $\bfPsi_x$, and the additional constraint (built into the definition of $\bfPsi_x^\dagger$) means that the unique shortest vector from the solution manifold is selected.

Hence, $\PR_{\fspace_N}{\KO}_{|\fspace_N} : \mathcal{F}_N \longrightarrow \mathcal{F}_N$, defined as 
$$
\PR_{\fspace_N}{\KO}_{|\fspace_N} \bfun_i = \sum_{j=1}^N \bfu_{ji}\bfun_j ,\;\; i=1, \ldots, N,
$$
is in the basis $\mathcal{B}=\{ \bfun_1, \ldots, \bfun_N\}$ represented by the uniquelly determined matrix
\begin{eqnarray}\label{eq:UN}
\!\!\!\!\!\!&&\!\! [\PR_{\fspace_N}\KO_{|\mathcal{F}_N}]_{\mathcal{B}} = \bfPsi_x^{\dagger}\bfPsi_y \equiv U_N = (\bfu_{ji})_{j,i=1}^N,\\
\!\!\!\!\!\!&&\!\!  \PR_{\fspace_N}\KO_{|\mathcal{F}_N}\! (\begin{pmatrix} \bfun_1(x) & \ldots & \bfun_N(x)\end{pmatrix}\!\! \begin{pmatrix}\mathsf{f}_1 \cr \vdots\cr \mathsf{f}_N\end{pmatrix})\! =\! \begin{pmatrix} \bfun_1(x) & \ldots & \bfun_N(x)\end{pmatrix} (U_N\!\! \begin{pmatrix}\mathsf{f}_1 \cr \vdots\cr \mathsf{f}_N\end{pmatrix}).\nonumber
\end{eqnarray}
If we introduce the notation $\bfunv (x)=( \bfun_1(x),\ldots, \bfun_N(x))^T\in\Cz^N$, 
and similarly for the residuals, $\bfrho(x)=( \rho_1(x), \ldots, \rho_N(x))^T\in\Cz^N$,
then the action of $\KO$ on $\fspace_N$ can be compactly written as
\begin{equation}
    \KO (\bfunv(x)^T \bfz) = \bfunv(x)^T U_N \bfz + \bfrho(x)^T \bfz,\;\;x\in\Cz^n,\; \bfz\in\Cz^N.
\end{equation}
In the sequel, we ease the notation and write
$\KO \bfunv(x)^T \bfz$ instead of $\KO (\bfunv(x)^T \bfz)$, and for a sequence of functions $\zeta_1(x), \zeta_2(x), \ldots$ we define
$\KO (\zeta_1(x), \zeta_2(x), \ldots)$ as
$(\KO\zeta_1(x), \KO\zeta_2(x), \ldots)$.
%

%
{
\subsection{Convergence as $m\rightarrow\infty$}\label{SS=Convergence}
In this section we discuss the convergence of $\PR_{\fspace_N}\restr{\KO}{\fspace_N} : \mathcal{F}_N \longrightarrow \mathcal{F}_N$ to ${\KO}$ as $N\rightarrow \infty$. As shown in \cite{korda-mezic-EDMD-convergence} it is enough to show this convergence in order to establish convergence of the EDMD operator $\PR_{\fspace_{N},m}\restr{\KO}{\fspace_N} : \mathcal{F}_N \longrightarrow \mathcal{F}_N$  to $\KO$ as $N,m\rightarrow \infty$\footnote{The discussion of this issue in \cite{korda-mezic-EDMD-convergence} pertains to the functional space being $L^2(\mu)$, but can be established for any separable Hilbert spaces adapted to dynamical systems, considered in \cite{mezic:2020}.}. 

For an operator $A$, we denote by $R_A(\lambda)=(A-\lambda I)^{-1}, \lambda \in \rho(A), R_A:\rho(A)\rightarrow \cB(\cH)$, its resolvent function, where $\rho(A)=\mathbb{C}/\sigma(A)$ is the resolvent set, $\sigma(A)$ is the spectrum of $A$, and $\cB(\cH)$ the space of bounded linear operators on $\cH$. We start with the notion of strong resolvent convergence:
\begin{definition}
A sequence of operators $B_n \in \cB(\cH)$ converges strongly to $B \in \cB(\cH)$ provided $B_n h\rightarrow Bh,\ \  \forall h\in \cH.$
\end{definition}
We assume $\PR_{\fspace_N}\rightarrow I$ strongly with $N$, where $I$ is the identity operator on $\cH$ (see \cite{korda-mezic-EDMD-convergence} for some conditions on the basis functions and $\cH$  under which this is true).
\begin{definition}
Let $A_n\in \cB(\cH), n\in \mathbb{N}$ and $A \in \cB(\cH)$. We say that $A_n$ converges to $A$ in the strong resolvent sense if $R_{A_n}\rightarrow R_{A} $ strongly for any - and thus every - $\lambda \in \Gamma$ where 
\begin{equation}
  \Gamma =\mathbb{C}/(\sigma(A)\cup_{n}\sigma(A_n))=\mathbb{C}/\Sigma.
\end{equation}
\end{definition}

\begin{lemma}
Assume $\KO$ admits the Kato spectral expansion
\begin{equation}
    \KO=\sum_{j}\lambda_jP^s_j+D_j
\end{equation}
where $D_j$ is quasinilpotent.
Then
$\PR_{\fspace_N}\restr{\KO}{\fspace_N}$ converges with $N$ to ${\KO}$  in the strong resolvent sense.
\end{lemma}
{\sc Proof:} The projections 
$\KO_N$ admit the expansion
\begin{equation}
    \KO_N=\sum_{j}\lambda_{j}P_NP^s_{j}P_N+P_ND_jP_N
    \label{eq:expN}
\end{equation}
where each eigenvalue $\lambda_j$ has finite algebraic multiplicity $m_j$ and the resolvent operator is
$$R_{\KO_N}(\lambda)=\sum_{j}\frac{P_NP^s_{j}P_N}{\lambda-\lambda_j}+\sum_{n=1}^{m_j-1}\frac{P_ND^n_{j}P_N}{(\lambda-\lambda_j)^{-n-1}}.$$
 Thus, $R_{\KO_N}$ strongly converges to $R_{\KO}(\lambda)$.\hfill$\boxtimes$

\begin{proposition}
Assume the spectrum of $\KO$ is contained inside the unit circle. Let $\KO_N$ converge with $N$ to ${\KO}$  in the strong resolvent sense. then $f(\KO_N)$ converges
to $f(\KO)$ strongly for every bounded $C^\infty$ function $f :\Sigma\rightarrow \mathbb{C}$.
\end{proposition}
{\sc Proof:} We need the following Lemma:
\begin{lemma} The spectrum of $\KO_N$ is inside the unit circle.
\label{l}
\end{lemma}
{\sc Proof:} Follows from   (\ref{eq:expN}) .\hfill$\boxtimes$

Now we finish the proof of the Proposition: There is a $C^\infty$ functional calculus for scalar operators \cite{maeda1962function}, that, together with stability established in the Lemma \ref{l}, implies the result. \hfill$\boxtimes$

Now we can discuss spectral convergence issues. We have the following:
\begin{theorem}
The sequence $\KO_N$ is spectraly inclusive, i.e. $\sigma (\KO)\subseteq \lim_{N\rightarrow \infty} \sigma (\KO_N)$.
\end{theorem}
{\sc Proof:} Suppose not. Then there is a $\lambda \in \mathbb{C}$ and $\epsilon>0$ such that $\sigma(\KO_N)\cap B_\epsilon(\lambda)=\emptyset$  where $B_\epsilon(\lambda) \subset \mathbb{C}$ is an open ball of radius $\epsilon$ centered at $\lambda$. Consider a $C^\infty$ function $f$ that is 1 on $B_{\epsilon/2}(\lambda)$ and zero outside of $B_\epsilon(\lambda)$. Then $f(\KO_N)=0$ implying $f(\KO)g=\lim f(\KO_N)g=0$ for every $g$. Since $\lambda \in \sigma (\KO)$ we have that the range of the spectral projection $P(B_{\epsilon/2}(\lambda))$ is not empty, containing a function $h$ such that $f(\KO)h=h$, leading to a contradiction.
\hfill$\boxtimes$
}

\subsection{The Rayleigh quotient of $U_N$}\label{SS=RQ_UN}

The  $N\times N$ matrix $U_N=\bfPsi_x^{\dagger}\bfPsi_y$ is in in Section \ref{SS=Num-Compress} defined as a matrix representation of the compression  $\PR_{\fspace_N}{\KO}_{|\fspace_N} : \mathcal{F}_N \longrightarrow \mathcal{F}_N$, where the projection is executed using algebraic least squares -- in the limited data driven setting this is a reasonable approximation. But, if   $N >  m$, it is of low rank of at most $m$, and potentially of large dimension. Moreover, in this case $\bfPsi_x$ and $\bfPsi_y$ are also of rank of at most $m$, and the data does not contain information to uniquely define $U_N$ in the sense of the formulation in Section \ref{SS=Num-Compress}. 

Indeed, if $E$ is such that $\bfPsi_x E=\0$, then
$\bfPsi_x (U_N + E)= \bfPsi_x U_N = \bfPsi_x\bfPsi_x^\dagger\bfPsi_y = \PR_x \bfPsi_y$, where $\PR_x$ is the orthogonal projector onto the range of $\bfPsi_x$.
The formula
$U_N=\bfPsi_x^{\dagger}\bfPsi_y$ is just a particular way to uniquely select a point in the linear manifold of the solutions of (\ref{eq:LS-Ui}, \ref{eq:LS-Ui-2}). It is appealing because 
of the minimal norm property and the corresponding theory of the Moore-Penrose generalized inverse, but it is not clear why should this choice be considered intrinsic.

The key observation is that the ambiguity in defining the matrix of the compressed operator is removed by taking the Rayleigh quotient of $U_N$ with respect to the range of $\bfPsi_x$:
$$
(\bfPsi_x^\dagger)^\dagger ( U_N + E) \bfPsi_x^\dagger = \bfPsi_x U_N \bfPsi_x^\dagger = \bfPsi_x\bfPsi_x^\dagger\bfPsi_y \bfPsi_x^\dagger = \PR_x \bfPsi_y \bfPsi_x^\dagger .
$$
For the reasons of numerical stability,  the matrix representation of the Rayleigh quotient should be computed using a unitary basis of  the range of $\bfPsi_x^\dagger$, or of some of its subspaces. Further, even when $m>N$ and all matrices are of full rank $N$, $\bfPsi_x$ can be ill-conditioned so that the numerical solution of the least squares problem becomes difficult and the numerical rank must be carefully determined.

To find a unitary basis and determine the correct dimension, a unitary rank revealing decomposition is used. The SVD is preferred because of its theoretical and numerical properties, and because of availability of robust and computationally efficient software implementations. 
Let $r=\mathrm{rank}(\bfPsi_x)$ and let the SVD of $\bfPsi_x$ be 
\begin{equation}\label{eq:SVD-psix}
\bfPsi_x = W \Sigma V^*,\;\;\Sigma=\mathrm{diag}(\sigma_i)_{i=1}^r,\, W\in\Cz^{m\times r},\;
V\in\Cz^{N\times r},\;W^*W=V^*V=\Id_r.
\end{equation}
(Here we write the \emph{economy size} SVD, i.e. the zero singular values and the corresponding singular vectors are chopped out from the full SVD.)

Using $\bfPsi_x^\dagger = V\Sigma^{-1}W^*$, the Rayleigh quotient $\widehat{U} = V^* U_N V$ is the $r\times r$ matrix
\begin{equation}\label{eq:RQ0}
\widehat{U} = V^* U_N V = V^*(V\Sigma^{-1}W^*)\bfPsi_y V = \Sigma^{-1}W^*\bfPsi_y V.
\end{equation}
It can be easily checked that $U_N V=V\widehat{U}$. 
We can express $V$ as $V=\bfPsi_x^* W\Sigma^{-1}$ and then $\widehat{U}$ as
\begin{equation}\label{eq:RQ}
\widehat{U} = V^* U_N V = V^*(V\Sigma^{-1}W^*)\bfPsi_y (\bfPsi_x^* W\Sigma^{-1}) = \Sigma^{-1}W^*(\bfPsi_y \bfPsi_x^*) W\Sigma^{-1}.
\end{equation}
The reason for considering (\ref{eq:RQ}) instead of the simpler formula $\widehat{U} = V^* U_N V$  will become apparent in Section \ref{SS=EDMD-kernek-trick}.

\subsubsection{Why is the SVD truncated?}
In practical finite precision computation, the rank of a matrix is determined as the numerical rank, based on the computed singular values and guided by the Eckart-Young-Mirsky theorem. Since this involves truncation of small (non-zero) singular values, the numerical rank could differ from the true rank. It should be stressed that determination of numerical rank is a nontrivial problem, see e.g. 
\cite{Golub:1976:RDL:892104}.

In the context of the Rayleigh-Ritz procedure of computing the eigenvalues and eigenvectors based on the
available data $\bfPsi_x$, $\bfPsi_y$, the issue is not about getting a good low rank approximation 
of $\bfPsi_x$ but about extracting spectral information based on the pair $\bfPsi_x$, $\bfPsi_y$. The SVD is used to
compute an orthonormal basis so that the Rayleigh-Ritz algorithm can be numerically more stable.
Using the low-rank approximation means discarding information that is in the range of $\bfPsi_x$, so it is better to get as large as possible orthonormal basis, i.e. more left singular vectors. 
Unfortunately, the smallest singular values (and the corresponding singular vectors) are often computed very inaccurately and the data driven formula (\ref{eq:RQ0}) for the Rayleigh quotient is contaminated by large errors if small singular values are used.
If a more accurate SVD algorithm capable of getting even small singular values accurately (when feasible) is available, then they can be included. Another important detail is that small singular values are a consequence of 
unfavourable scaling of the data matrices. These important issues are discussed in detail in 
\cite{DDMD-SISC-2018}, \cite{Drmac-2020-koopman-book-chapter},   \cite{drmac-DMD-2023}, \cite{drmac-SYDMD-2023},  (For software details see \cite{lawn-298-drmac}, \cite{lawn-300-drmac}.)

Another important issue is the possibility of losing (through truncation) low energy unstable modes, see e.g. a discussion in \cite{ahuja_rowley_2010}.

\subsection{The EDMD with the kernel trick}\label{SS=EDMD-kernek-trick}

The key issue is practical computation when the data snapshots are mapped in ever higher dimensional space, i.e. when $N\rightarrow\infty$. The matrices $W$ and $\Sigma$ in (\ref{eq:SVD-psix}), (\ref{eq:RQ}) can be computed from the spectral decomposition of the  $m\times m$ cross-product  $C_{xx}=\bfPsi_x\bfPsi_x^*$, since $C_{xx} = W \Sigma^2 W^*$. The cross-product matrix $C_{yx}=\bfPsi_y\bfPsi_x^*$, used in (\ref{eq:RQ}), is $m\times m$ and its elements, as those of $C_{xx}$, are inner products of vectors from $N$-dimensional space. Hence, the only place where the potentially extremely large dimension $N$ incurs computational cost is in computing $C_{xx}$ and $C_{yx}$. 

Note that $C_{xx}=C_{xx}^*$, and that for all $i, j$ the elements of these matrices can be expressed using the Euclidean inner product $\langle\cdot,\cdot\rangle_{\Cz^N}$ as
\begin{eqnarray}
(C_{xx})_{ij} &=& (\bfPsi_x\bfPsi_x^*)_{ij} = \sum_{k=1}^N \bfun_k(\x_i)\overline{ \bfun_k(\x_j)} = \langle\bfunv(\x_i),\bfunv(\x_j)\rangle_{\Cz^N} , \\ (C_{yx})_{ij} &=& (\bfPsi_y\bfPsi_x^*)_{ij} = \sum_{k=1}^N \bfun_k(\y_i)\overline{ \bfun_k(\x_j)}  = \langle\bfunv(\y_i),\bfunv(\x_j)\rangle_{\Cz^N} .
\end{eqnarray}
This is where the kernel trick comes to the rescue. The trick is to work with the basis functions $\bfun_k$ implicitly by using a suitable function $\kerf$ such that
\begin{equation}
\kerf(x,y) = \langle \bfunv(x), \bfunv(y)\rangle_{\Cz^N},\;\;x, y\in\Cz^n ,   
\end{equation}
so that $(C_{xx})_{ij}=\kerf(\x_i,\x_j)$, $(C_{yx})_{ij}=\kerf(\y_i,\x_j)$. Since the construction starts with the function $\kerf$, its defining property must be the definiteness in the sense of the Mercer's theorem.
With a broad palette of kernel functions, so that computation of the Rayleigh quotient (\ref{eq:RQ}) is feasible for large dimensions, even when $N\rightarrow\infty$. For more details we refer to \cite{williams-2015-EDMD}.
%
%
%

\begin{remark}
{\em
The kernel trick is an excellent example of trade-off between computational efficiency and numerical accuracy. The key of the accuracy of the computed SVD
of $\bfPsi_x$ is the condition number $\kappa_2(\bfPsi_x)=\|\bfPsi_x\|_2 \|\bfPsi_x^\dagger\|_2$.
On the other hand, the condition number of $C_{xx}=\bfPsi_x\bfPsi_x^*$ is squared,  $\kappa_2(C_{xx})=\kappa_2(\bfPsi_x)^2$.
}
\hfill $\boxtimes$\end{remark}

\subsubsection{Connection with the DMD}\label{SS=DMD-connection-1}
The DMD is a special case of EDMD with $N=n$ and with the coordinate functions $\bfun_i(x)=e_i^T x$ as basis functions. The kernel function is $\kerf(x,y)=\langle x, y\rangle_{\Cz^n}$.

In the DMD, the snapshots are collected column-wise into the matrices $\X=(\x_1\;\ldots\; \x_m)$, $\Y=(\y_1\;\ldots\; \y_m)$. Note that $\X=\bfPsi_x^T$, $\Y=\bfPsi_y^T$. 
The DMD matrix is $A_n=\Y\X^\dagger=\bfPsi_y^T\bfPsi_x^{T\dagger}$, hence $A_n=A_N=U_N^T$. 
The SVD of $\X$ (using (\ref{eq:SVD-psix})) reads $\X = V^{*T}\Sigma W^T$, so that the $r\times r$ Rayleigh quotient $\widehat{A}$  used in the DMD is
$$
\widehat{A}=V^T A_N V^{*T} = (V^* U_N V)^T = \widehat{U}^T .
$$
Here too, in numerical computations, $r$ is the numerical rank, and $\widehat{A}$ is computed using the formula
\begin{equation}
\widehat{A}=V^T A_N V^{*T} = V^T \Y W^{T*}\Sigma^{-1} V^T V^{*T} = V^T \Y W^{T*}\Sigma^{-1} .
\end{equation}
In DMD, a spectral decomposition of $\widehat{A}$, $\widehat{A}= G \Lambda G^{-1}$, is computed and the approximate eigenvectors of $A_N$ are
the columns of $V^{*T} G$.

For a numerically robust implementation of the DMD see \cite{DDMD-SISC-2018}, \cite{lawn-298-drmac} and for a special case of Hermitian DMD see \cite{lawn-300-drmac}. 
\begin{remark}
{\em 
The DMD is often in the literature interpreted as a mere linear regression problem, i.e. 
$A_n=\Y\X^\dagger$ is taken as a mere solution of the linear least squares problem $\| \Y - A_n \X\|_F\rightarrow\min$, 
with the solution $A_n=\Y\X^\dagger$, and ignoring the connection with the Koopman linearization of the underlying dynamics and
the fact that $A_n$ is, in the case of tall data matrices, just one out of many solutions of the LS problem. 
See \cite{drmac-DMD-2023}, \cite{drmac-SYDMD-2023} for more detailed discussion.\footnote{See also  \cite{lawn-298-drmac}, \cite{lawn-300-drmac}.}
}
\hfill $\boxtimes$\end{remark}
\subsection{Evaluation of the approximate eigenfunctions of $\KO$}\label{SS=evaluate-fs}
If we compute the $r$ eigenpairs $(\omega_i,s_i)$ of $\widehat{U}$, $\widehat{U}s_i = \omega_i s_i$, $i=1,\ldots, r$, then 
$$
U_N Vs_i = \bfPsi_x^\dagger \bfPsi_y V s_i = V \Sigma^{-1}W^*\bfPsi_y\bfPsi_x^* W\Sigma^{-1} s_i = \omega_i V s_i ,
$$
so that $(\omega_i, Vs_i)$ are the eigenpairs of $U_N$. 
The approximate eigenfunctions of $\KO$ are
\begin{equation}\label{eq:eigf_U}
\eigf_i(x) = \begin{pmatrix}
\bfun_1(x) & \bfun_2(x) & \ldots& \bfun_N(x)
\end{pmatrix} Vs_i  = \sum_{j=1}^N \bfun_j(x) (Vs_i)_j,\;\;i=1,\ldots, r.
\end{equation}
This is understood in the sense that
\begin{eqnarray}
\KO \eigf_i(x) &=&   \begin{pmatrix}
\bfun_1(x) & \bfun_2(x) & \ldots& \bfun_N(x)
\end{pmatrix} U_N  V s_i  + \bfrho(x)^T Vs_i\\
&=& \begin{pmatrix}
\bfun_1(x) & \bfun_2(x) & \ldots& \bfun_N(x)
\end{pmatrix} \omega_i V s_i + \bfrho(x)^T Vs_i\\  &=& \omega_i \eigf_i(x) + \bfrho(x)^T Vs_i \backsimeq \omega_i \eigf_i(x).
\end{eqnarray}
(Recall that $\backsimeq$ is a shorthand that means equal up to a residual that is minimized over the available data in the sense of least squares.)

In a practical computation, the $\eigf_i$'s are not numerically evaluated using the definition (\ref{eq:eigf_U}), as it involves the dimension $N$.

\subsubsection{Evaluation at the snapshots $\x_j$}
Evaluating the eigenfunctions $\eigf_i$ at the snapshots $\x_j$ means tabulating the values
$$
\bfPhi_x = \begin{pmatrix}
\eigf_1(\x_1) & \ldots & \eigf_r(\x_1) \cr
\vdots & \ddots & \vdots \cr\eigf_1(\x_m) & \ldots & \eigf_r(\x_m)
\end{pmatrix} = \bfPsi_x V S \in\Cz^{m\times r} ,
$$
where the relation $\bfPhi_x = \bfPsi_x V S$ follows directly from the definition of the matrix $\bfPsi_x$ and the eigenfunctions (\ref{eq:eigf_U}). Next, observe that the SVD $\bfPsi_x=W\Sigma V^*$ eliminates $V$ so that we have $\bfPhi_x = W\Sigma S$. In practice, the SVD of $\bfPsi_x$  may not be feasible and the kernel trick is used, so that the same conclusion follows from the spectral decomposition of $\bfPsi_x \bfPsi_x^*$ and the relation $V=\bfPsi_x^* W\Sigma^{-1}$ :
$$
\bfPhi_x = \bfPsi_x \bfPsi_x^* W \Sigma^{-1} S = W\Sigma^2 W^* W\Sigma^{-1} S = W \Sigma S.
$$
\subsubsection{Evaluation at an arbitrary $\x\in\Cz^n$}\label{SS=phi(any_x)}
To evaluate $\eigf_i(\x)$ at an arbitrary $\x\in\Cz^n$ we need
$$
\eigf_i(\x) = \bfunv(\x)^T V s_i = \bfunv(\x)^T \bfPsi_x^* W\Sigma^{-1} s_i ,
$$
where, using the kernel function $\kerf(\cdot,\cdot)$, 
$$
(\bfunv(\x)^T \bfPsi_x^*)_j = 
\sum_{k=1}^N \bfun_k(\x)\overline{\bfun_k(\x_j)} = \langle \bfunv(\x), \bfunv(\x_j)\rangle_{\Cz^N} = \kerf (\x,\x_j) .
$$
Hence
$$
\eigf_i(\x) = \begin{pmatrix} 
\kerf(\x,\x_1) & \kerf(\x,\x_2) & \ldots & \kerf(\x,\x_m)
\end{pmatrix} W\Sigma^{-1} s_i .
$$
%
%
\subsection{Snapshot representation}\label{SS=EDMD-snap-repr}
Next computational task is to represent the snapshots $\x_1, \ldots, \x_m$ using the computed eigenfunctions $\eigf_1, \ldots, \eigf_r$.
More generally, given $\ell$ scalar functions $g_1, \ldots, g_\ell$, we seek an approximation
$$
\begin{pmatrix} g_1(x) & \ldots & g_\ell(x)\end{pmatrix} \backsimeq \begin{pmatrix}\eigf_1(x) & \ldots & \eigf_r(x) \end{pmatrix} \Gamma,\;\;
\Gamma\in\Cz^{r\times\ell}.
$$
The approximation is computed in the least squares sense, column-wise, e.g.
\begin{equation}\label{eq:Gamma_i}
\Gamma(:,i) = \mathrm{argmin}_{h\in\Cz^r} \| \left(\begin{smallmatrix} g_i(\x_1)\cr \vdots\cr g_i(\x_m)\end{smallmatrix}\right) - \bfPhi_x h\|_2 = \bfPhi_x^\dagger \left(\begin{smallmatrix} g_i(\x_1)\cr \vdots\cr g_i(\x_m)\end{smallmatrix}\right) ,\;\;i=1,\ldots, \ell.
\end{equation}
For snapshot reconstruction, $\ell=n$ and $g_i(x)=e_i^T x$ is the $i$-th coordinate function, and in the sequel we consider this special case. If we define $\X=(\x_1\; \ldots \; \x_m)$ then the column-wise definition (\ref{eq:Gamma_i}) of $\Gamma$ gives (in the least squares sense) $\X^T \approx \bfPhi_x\Gamma$ and the solution matrix $\Gamma$ can be compactly written as
$$
\Gamma = \bfPhi_x^\dagger \X^T = (W\Sigma S)^\dagger \X^T = S^{-1}\Sigma^{-1}W^* \X^T .
$$
Hence, we have the representation of the snapshots (in a least squares sense)
$$
\X = (\x_1\; \ldots \; \x_m) \approx  \Gamma^T\bfPhi_x^T,\;\;\Gamma^T = \X W^{*T}\Sigma^{-1} S^{-T}= ( \gamma_1\; \ldots \; \gamma_r).
$$
Hence, each snapshot has a modal decomposition
$$
\x_i \approx \sum_{j=1}^r \gamma_j \eigf_j(\x_i),\;\;i=1,\ldots, m.
$$
Note that $\widehat{U} = S \Omega S^{-1}$ yields 
$\widehat{U}^T = S^{-T} \Omega S^{T}$, i.e. the  columns of $S^{-T}$ are the eigenvectors of $\widehat{U}^T$. (In the case of DMD, $\widehat{U}^T=\widehat{A}$; see Section \ref{SS=DMD-connection-1}.)

\begin{remark}\label{REM-caveat}
{\em 
Note that both $S$ and its inverse appear in the formuals -- for the values $\eigf_k(\x_i)$ (the matrix $\bfPhi_x$) we need $S$, and for the modes $\gamma_j$  (the matrix $\Gamma$) we need $S^{-T}$. 
There is an important caveat here:  $S$ may be severely ill-conditioned (even numerically singular) and the method requires using $S^{-1}$. This important issue is addressed in Section \ref{S=KSD}.
}
\hfill $\boxtimes$\end{remark}

\subsubsection{Connection to the Schmid's DMD}

In the case of DMD (recall Section \ref{SS=DMD-connection-1}), $\bfPsi_x = \X^T = W\Sigma V^*$, $\X = V^{*T}\Sigma W^T$ and 
$$
\Gamma^T =  V^{*T}\Sigma W^T W^{*T}\Sigma^{-1} S^{-T} = V^{*T} S^{-T} = V^{*T} G .
$$
\section{Koopman-Schur decomposition}\label{S=KSD}
In the EDMD and the DMD (see Sections \ref{SS=DMD-connection-1}, \ref{SS=evaluate-fs}) computational frameworks it is tacitly assumed that the Rayleigh quotient matrix $\widehat{U}$ ($\widehat{A}$ in the case of DMD) is diagonalizable, and the whole construction is build on the decomposition 
\begin{equation}\label{eq:U=SOS-1}
\widehat{U} = S \Omega S^{-1},\;\;\Omega=\mathrm{diag}(\omega_i)_{i=1}^r,\;\; S = (s_1,\ldots,s_r),\;\widehat{U} s_i = \omega_i s_i. 
\end{equation}
Here $\widehat{U}=V^* U_N V$, where $U_N$ represents a compression of $\KO$ onto the range of the orthonormal matrix $V$; see (\ref{eq:RQ}) in  Section \ref{SS=RQ_UN}.

In Section \ref{S=Intro}, we used a $2\times 2$ example to illustrate the problem of non-existence of the full system of eigenvectors. Furthermore, in practical computation, although diagonalizability of a matrix is not a priori guaranteed, a numerical algorithm will most probably return the maximal number of $r$ eigenvectors, even if there exists only one. The following example will illustrate.
\begin{example}
{\em 
We use the Matlab's function \texttt{eig} to compute the eigenvalues and eigenvectors of a Jordan block of dimension $10$. The computed eigenvector matrix $S$ will be singular, without any warning/error message, as shown by the following simple code
\begin{verbatim}
>> Uhat=diag(ones(10,1),1);
>> [S,~]=eig(Uhat); [G,~]=eig(Uhat');
>> [ cond(S) cond(G) ]
ans =
   Inf   Inf 
\end{verbatim} 
More generally, we can think of $\widehat{U}$ being equal or close to a matrix whose Jordan normal form contains several (Jordan) blocks of different sizes, corresponding to same and/or different eigenvalues.
 As shown in Sections \ref{SS=DMD-connection-1} and \ref{SS=EDMD-snap-repr}, the eigenvector matrices $S$ and $G$ are the ingredients of the matrix of the Koopman modes, which is in this case of numerical rank one (in the general case possibly highly ill-conditioned and numerically rank deficient) and as such it has a rather limited potential for discovering latent properties.
 }
\end{example}
Recall that matrices with non-diagonal Jordan form are a nowhere dense set, so it is rather unlikely that our $\widehat{U}$ is such. However, even if diagonalizable,  $\widehat{U}$ can be highly non-normal and the matrix $S$ can be numerically ill-conditioned, and the eigenvectors can  be so volatile that using them numerically may not have sense at all. On the other hand, those eigenvectors are the key ingredients in the modal decomposition in the EDMD -- recall Remark \ref{REM-caveat}.

It is also possible that some individual eigenvectors themselves are sensitive to perturbations, but a spectral subspace that encompasses them is much more robust. Multiple or closely clustered eigenvalues that make a cluster well separated from the rest of the spectrum have such eigenvectors.  %

These nontrivial numerical difficulties have been the main motivation and a driving force for the development of an alternative approach that is based on the ordered partial Schur decomposition.
The key idea is to avoid the spectral decomposition (diagonalization, Jordan normal form) of the compressed operator $\KO$ (see Section \ref{SS=Num-Compress}) and to use numerically more robust Schur decomposition. From a numerical linear algebra perspective, this seems plausible and it only remains to adapt this approach to the DMD/EDMD and to ensure the same functionality in the applications. 
In this section we provide the details of the new approach. 

\subsection{Partial Schur decomposition of $\KO$ and $U_N$}
Following the above discussion, let us, instead of the diagonalization (\ref{eq:U=SOS-1}), compute the Schur form  (Theorem \ref{TM:Schur-form}) of $\widehat{U}$:
\begin{equation}\label{eq:Schur-hat-U}
    \widehat{U} = Q T Q^*, \;\;T=\left(\begin{smallmatrix}
    t_{11} & \ldots & t_{1r} \cr
           & \ddots  & \vdots \cr
           &        & t_{rr}
    \end{smallmatrix}\right),\;\;t_{ii}=\omega_i,
\end{equation}
where $Q$ is unitary and $T$ is upper triangular with the eigenvalues of $\widehat{U}$, in some order, on the diagonal. 
Then
\begin{equation}\label{eq:VQ-partial-Schur}
    U_N (VQ) = \bfPsi_x^\dagger\bfPsi_y VQ = V\Sigma^{-1}W^* \bfPsi_y\bfPsi_x^*W\Sigma^{-1} Q = V\widehat{U} Q = (V Q) T
\end{equation}
represents a partial Schur form  of $U_N$ -- the matrix $VQ$ is $N\times r$ with $N  > r$ (including the case of $N=\infty \gg r$). Its range is an $r$-dimensional $U_N$-invariant subspace and the corresponding compression of $U_N$ is $T$. 
If we set $Z=VQ$, then 
\begin{equation}\label{eq:UZ=ZT}
U_N Z = ZT,\;\;\mbox{and, more generally,}\;\;U_N^k Z = Z T^k,\;\;k=1, 2, \ldots
\end{equation}
In fact, $f(U_N)Z = Z f(T)$, where $f$ is any analytic function defined in a neighborhood of the spectrum of $U_N$. Note that $Z^* Z = Q^* V^* V Q = \Id_r$.

The relation $U_N Z=ZT$ can be illustrated as 
$$
\left(\begin{smallmatrix}
\square & \square & \square & \square\cr
\square & \square & \square & \square\cr
\square & \square & \square & \square\cr
\square & \square & \square & \square\cr
\end{smallmatrix}\right)
\left(\begin{smallmatrix}
\blacksquare & \blacksquare \cr
\blacksquare & \blacksquare \cr
\blacksquare & \blacksquare \cr
\blacksquare & \blacksquare \cr
\end{smallmatrix}\right)
= \left(\begin{smallmatrix}
\blacksquare & \blacksquare \cr
\blacksquare & \blacksquare \cr
\blacksquare & \blacksquare \cr
\blacksquare & \blacksquare \cr
\end{smallmatrix}\right)
\left(\begin{smallmatrix}
\boxtimes & \boxtimes \cr
 \mathbf{0} & \boxtimes
\end{smallmatrix}\right) ,
$$
which explains the term \emph{partial Schur decomposition} of $U_N$ -- the triangular factor $T$ is of smaller dimension than $U_N$ and only a
portion of the eigenvalues is computed.
\subsubsection{Schur functions of $\KO$}
On the operator level, the relation (\ref{eq:UZ=ZT}) becomes 
\begin{eqnarray}
\!\!\!\!\!\KO (\begin{pmatrix}
\bfun_1(x) & \ldots & \bfun_N(x)
\end{pmatrix} Z)  &=&  \begin{pmatrix}
\bfun_1(x) & \ldots & \bfun_N(x)
\end{pmatrix} U_N Z +\bfrho(x)^T Z  \nonumber\\
&=& \begin{pmatrix}
\bfun_1(x) & \ldots & \bfun_N(x)
\end{pmatrix} Z T +\bfrho(x)^T Z \label{eq:residual-rhoTZ}\\
&\backsimeq& \begin{pmatrix}
\bfun_1(x) & \ldots & \bfun_N(x)
\end{pmatrix} Z T .
\end{eqnarray}
If we define a new sequence of functions as
\begin{equation}
    \begin{pmatrix}
    \zeta_1(x) & \zeta_2(x) & \ldots & \zeta_r(x)
    \end{pmatrix} = \begin{pmatrix}
\bfun_1(x) & \ldots & \bfun_N(x)
\end{pmatrix} Z ,
\end{equation}
then we have, using (\ref{eq:residual-rhoTZ}), 
\begin{eqnarray}
\!\!\!\!    \KO \begin{pmatrix}
    \zeta_1(x) & \zeta_2(x) & \ldots & \zeta_r(x)
    \end{pmatrix} &=& \begin{pmatrix}\zeta_1(x) & \zeta_2(x) & \ldots & \zeta_r(x)
    \end{pmatrix} T +\bfrho(x)^T Z \nonumber\\ 
    &\backsimeq& \begin{pmatrix}
    \zeta_1(x) & \zeta_2(x) & \ldots & \zeta_r(x)
    \end{pmatrix} T . \label{eq:Uzeta=zetaT(x)}
\end{eqnarray}
Since $T$ is upper triangular, relation (\ref{eq:Uzeta=zetaT(x)}) contains a nested sequence of analogous partial triangulations 
\begin{equation}\label{eq:U-Schur-i}
    \KO \begin{pmatrix}
    \zeta_1 & \zeta_2 & \ldots & \zeta_i
    \end{pmatrix} \backsimeq \begin{pmatrix}
    \zeta_1 & \zeta_2 & \ldots & \zeta_i
    \end{pmatrix} T(1:i,1:i) , \;\;i=1, \ldots, r.
\end{equation}
In particular, 
$\KO\; \mathrm{Span}\!\begin{pmatrix}
    \zeta_1 & \zeta_2& \ldots & \zeta_i
    \end{pmatrix} \precsim \mathrm{Span}\!\begin{pmatrix}
    \zeta_1 & \zeta_2& \ldots & \zeta_i
    \end{pmatrix} , \;\;i=1, \ldots, r,
$
where the relation $\precsim$ is used as a shorthand to be interpreted as follows: any function in the subspace on the left hand side can be written  as a linear combination of the $\zeta_i$'s plus a residual that has certain minimality property over the available data (See Section \ref{SS=Num-Compress}.).

\begin{remark}
{\em
 If the basis functions $\bfun_i$ are an orthonormal set,\footnote{If the basis function are not normalized (but are mutually orthogonal), we simply redefine them to $\bfun_i\leftarrow \bfun_i/\|\bfun_i\|_{L^2}$, to enforce $\langle \bfun_i,\bfun_j\rangle_{L^2}=\bfdelta_{ij}$.} $\langle \bfun_i,\bfun_j\rangle_{L^2}=\bfdelta_{ij}$, then 
\begin{equation}\label{eq:zeta-orth}
\langle \zeta_i,\zeta_j\rangle_{L^2} = \langle \sum_k (z_i)_k \bfun_k, \sum_\ell (z_j)_\ell \bfun_\ell\rangle_{L^2} = \sum_k (z_i)_k \overline{(z_j)_k} = \langle z_i, z_j\rangle = \bfdelta_{ij}.
\end{equation}
Otherwise, we can apply the Gram-Schmidt orthogonalization (in $L^2$)
$$
\begin{pmatrix}
    \widehat{\zeta}_1 & \widehat{\zeta}_2 & \ldots & \widehat{\zeta}_r
    \end{pmatrix} = 
\begin{pmatrix}
    \zeta_1 & \zeta_2 & \ldots & \zeta_r
    \end{pmatrix} R ,
$$
where  $\langle \widehat{\zeta}_i, \widehat{\zeta}_j\rangle_{L^2} = \bfdelta_{ij}$, and $R$ is upper triangular.
Then, $\widehat{T} = R^{-1} TR$ is upper triangular and we have
\begin{eqnarray}\label{eq:Uzeta=zetaT(x)-orth}
\!\!\!\!    \KO \begin{pmatrix}
    \widehat{\zeta}_1(x) & \widehat{\zeta}_2(x) & \ldots & \widehat{\zeta}_r(x)
    \end{pmatrix} &=& \begin{pmatrix}\widehat{\zeta}_1(x) & \widehat{\zeta}_2(x) & \ldots & \widehat{\zeta}_r(x)
    \end{pmatrix} \widehat{T} +\bfrho(x)^T VQ R \nonumber \\ 
    &\backsimeq& \begin{pmatrix}
    \widehat{\zeta}_1(x) & \widehat{\zeta}_2(x) & \ldots & \widehat{\zeta}_r(x)
    \end{pmatrix} \widehat{T} .
\end{eqnarray}
}
\hfill $\boxtimes$\end{remark}
%
%
\subsection{Ordered partial Schur form}\label{SS=ordered-Schur}
The partial Schur form can be truncated at any index $\ell<r$, resulting again in a partial Schur form of lower order, see (\ref{eq:U-Schur-i}). It may be of interest to have such a low order partial Schur form that contains only some selected eigenvalues that are not necessarily the leading $i$ values on the diagonal of $T$. 
On the other hand, when we compute the Schur form numerically, the order of the eigenvalues on the diagonal of $T$ is in general not known in advance -- it depends on the convergence of the numerical algorithm/software  for a given input $\widehat{U}$.

A given Schur form can be reordered in the sense that an efficient algorithm determines a unitary matrix $\Theta$ such that $\widetilde{T} = \Theta^* T \Theta$ is again upper triangular with diagonal entries corresponding to the eigenvalues in any given order; see e.g. \cite{BAI199375}.
The new partial Schur form of $U_N$ becomes
$$
U_N (VQ\Theta) = (VQ\Theta) \widetilde{T} ,\;\;\widetilde{T}=\Theta^* T \Theta, 
$$
and we replace (\ref{eq:Uzeta=zetaT(x)}) with 
\begin{equation}
    \KO \begin{pmatrix}
    \zeta_1(x) & \zeta_2(x) & \ldots & \zeta_r(x)
    \end{pmatrix}\Theta \backsimeq \begin{pmatrix}
    \zeta_1(x) & \zeta_2(x) & \ldots & \zeta_r(x)
    \end{pmatrix}\Theta (\Theta^* T \Theta) ,
\end{equation}
i.e. the new functions are generated using $\widetilde{Z}=Z\Theta = VQ\Theta$, 
\begin{equation}
\begin{pmatrix}
    \widetilde{\zeta}_1(x) &  \ldots & \widetilde{\zeta}_r(x)
    \end{pmatrix} = 
    \begin{pmatrix}
    \zeta_1(x) & \ldots & \zeta_r(x)
    \end{pmatrix}\Theta = \begin{pmatrix}
\bfun_1(x) & \ldots & \bfun_N(x)
\end{pmatrix} \widetilde{Z} .
\end{equation}
Since $\widetilde{Z}^*\widetilde{Z}=\Id_r$, the orthogonality similar to (\ref{eq:zeta-orth}) remains true. The new partial Schur decomposition $U_N \widetilde{Z}=\widetilde{Z}\widetilde{T}$ can be truncated as in (\ref{eq:U-Schur-i}) at some index $i$, with the desired eigenvalues on the diagonal of $\widetilde{T}(1:i,1:i)$, and the corresponding Schur functions $\widetilde{\zeta}_1,\ldots, \widetilde{\zeta}_i$.
\begin{remark}
{\em
It should be noted that the reordering of the Schur form works on the triangular factor that is of modest dimension, and computing the additional reordering transformation $\Theta$ is not prohibitively expensive.
LAPACK library provides functions \texttt{STRSEN, DTRSEN, CTRSEN, ZTRSEN} that implement Schur form reordering for real and complex matrices, single and double precision. Its application to $Q$ or $VQ$ is a BLAS 3 operation that can be executed very efficiently, or (depending on the application of the decomposition) the Schur vectors can be kept in factored form. Matlab package contains the function \texttt{ordschur} that implements a method by Kressner \cite{Kressner-ordschur-2006}.
}
\hfill $\boxtimes$\end{remark}
%
%
%
\subsection{Evaluation of the Schur functions}\label{SS:eval-Schur-f}
The kernel trick evaluation of the functions $\zeta_i(\cdot)$ follows the lines of Section \ref{SS=evaluate-fs}.
At an arbitrary $\x\in\Cz^n$ we have (see Section \ref{SS=phi(any_x)})
$\zeta_i(x)=\bfpsi(x)^T V Q e_i=\bfpsi(x)^T \bfPsi_x^* W\Sigma^{-1}Q e_i$
and thus
$$
\begin{pmatrix}\zeta_1(\x) & \ldots & \zeta_r(\x) \end{pmatrix} = \begin{pmatrix} 
\kerf(\x,\x_1) & \kerf(\x,\x_2) & \ldots & \kerf(\x,\x_m)
\end{pmatrix} W\Sigma^{-1} Q .
$$
Note that a kernel trick can be used to compute $\zeta_i(\x)$.
At the snapshots, we have the tabulated values computed as 

\begin{equation}\label{zetax-SVD}
\bfzeta_x = \begin{pmatrix}
\zeta_1(\x_1) & \ldots & \zeta_r(\x_1) \cr
\vdots & \ddots & \vdots \cr\zeta_1(\x_m) & \ldots & \zeta_r(\x_m)
\end{pmatrix} = \bfPsi_x V Q = W\Sigma^2 W^* W\Sigma^{-1} Q = W \Sigma Q \in\Cz^{m\times r} .
\end{equation}
Note that $\bfzeta_x$ is conveniently given factored in a SVD decomposition. 
%

\subsection{Representation of observables using Schur functions}\label{SS=Observ_Representation_Schur_f}
As in Section \ref{SS=EDMD-snap-repr}, we seek representation of given scalar functions $g_1, \ldots, g_\ell$, but this time in terms of  
the Schur functions $\zeta_1,\ldots, \zeta_r$,
\begin{equation}\label{eq:g-repr}
\begin{pmatrix} g_1(x) & \ldots & g_\ell(x)\end{pmatrix} \backsimeq \begin{pmatrix}\zeta_1(x) & \ldots & \zeta_r(x) \end{pmatrix} \Xi,\;\;
\Xi\in\Cz^{r\times\ell}.
\end{equation}
In the restrictive data driven framework, the representation (\ref{eq:g-repr}) is feasible only at the data snapshots $\x_j$, i.e. at the data matrix $G_{ji}=g_i(\x_j)$.
The approximation is computed in the least squares sense, column-wise, e.g.
$$
\Xi(:,i) = \mathrm{argmin}_{h\in\Cz^r} \| \left(\begin{smallmatrix} g_i(\x_1)\cr \vdots\cr g_i(\x_m)\end{smallmatrix}\right) - \bfzeta_x h\|_2 = \bfzeta_x^\dagger \left(\begin{smallmatrix} g_i(\x_1)\cr \vdots\cr g_i(\x_m)\end{smallmatrix}\right) ,\;\;i=1,\ldots, \ell.
$$
Hence, $G\approx \bfzeta_x \Xi$ with 
$$
\Xi = \bfzeta_x^\dagger G = \begin{pmatrix}
\zeta_1(\x_1) & \ldots & \zeta_r(\x_1) \cr
\vdots & \ddots & \vdots \cr\zeta_1(\x_m) & \ldots & \zeta_r(\x_m)
\end{pmatrix}^\dagger \begin{pmatrix}
g_1(\x_1) & \ldots & g_\ell(\x_1) \cr
\vdots & \ddots & \vdots \cr g_1(\x_m) & \ldots & g_\ell(\x_m)
\end{pmatrix} = Q^* \Sigma^{-1} W^* G,
$$
where the pseudoinverse $\bfzeta_x^\dagger$ is computed directly from the readily available SVD (\ref{zetax-SVD}) of $\bfzeta_x$.

In  a DMD framework, the data snapshots are usually stored column-wise in a matrix, so that the column index represents the discrete time. For such a representation, write $G^T\approx \Xi^T \bfzeta_x^T$, introduce column partition $\Xi^T = \begin{pmatrix} \eta_1 & \ldots & \eta_r\end{pmatrix}$ and write the representation (\ref{eq:g-repr}) at the snapshots $\x_i$ as 
\begin{equation}\label{eq:GT-reconstruct}
\begin{pmatrix}
 g_1(\x_1) & g_1(\x_2) & \ldots & g_1(\x_m)\cr
 g_2(\x_1) & g_2(\x_2) & \ldots & g_2(\x_m)\cr
 \vdots & \vdots & \cdots & \vdots \cr
 g_\ell(\x_1) & g_\ell(\x_2) & \ldots & g_\ell(\x_m)\cr
\end{pmatrix}
\approx  
\begin{pmatrix} \eta_1 & \ldots & \eta_r\end{pmatrix}\!\!
\begin{pmatrix}
 \zeta_1(\x_1) & \zeta_1(\x_2) & \ldots & \zeta_1(\x_m)\cr
 \zeta_2(\x_1) & \zeta_2(\x_2) & \ldots & \zeta_2(\x_m)\cr
 \vdots & \vdots & \cdots & \vdots \cr
 \zeta_r(\x_1) & \zeta_r(\x_2) & \ldots & \zeta_r(\x_m)\cr
\end{pmatrix} .
\end{equation}
The special case $G^T=\X$ (the case of full state observables) is described in the next two propositions.
\begin{proposition}\label{PROP:3.1}
For the data snapshots representation, we specify the observables $g_i$ as: $\ell=n$, $g_i(x)=\ecv_i^T x$. Then  $G=\X^T$, and we have 
$$
\Xi = \bfzeta_x^\dagger \X^T = (W\Sigma Q)^\dagger \X^T = Q^*\Sigma^{-1}W^* \X^T ,
$$
i.e.
\begin{equation}\label{eq:xi_approx}
\X \approx \Xi^T \bfzeta_x^T,\;\;\mbox{where}\;\;\Xi^T = \X W^{*T}\Sigma^{-1}Q^{*T},\;\;\bfzeta_x^T = Q^T\Sigma W^T,
\end{equation}
and, with the column partition $\Xi^T = \begin{pmatrix} \eta_1 & \ldots & \eta_r\end{pmatrix}$, 
$$
\x_i \approx \sum_{j=1}^r\eta_j \zeta_j(\x_i) .
$$
\end{proposition}
This representation of the data snapshot matrix $\X$ simplifies if the basis functions are $\bfun_i(x)=\ecv_i^T x$, i.e. when $\bfPsi_x=\X^T$.
\begin{proposition}\label{PROP:snap-rep-2}
In the case $\bfPsi_x=\X^T = W\Sigma V^*$ and $\kerf(x,y)=\langle x,y\rangle_{\Cz^n}$, the representation (\ref{eq:xi_approx}) of the data snapshot matrix $\X$ in Proposition \ref{PROP:3.1} holds with 
$$
\Xi^T = \X W^{*T}\Sigma^{-1}Q^{*T} = V^{*T} Q^{*T}= Z^{*T}=\begin{pmatrix} \eta_1 & \ldots & \eta_r\end{pmatrix} .
$$
Note that in this case $\Xi^T=Z^{*T}$ is orthonormal, where $Z=VQ$ is the orthonormal matrix from the partial Schur decomposition (\ref{eq:VQ-partial-Schur}). 
Further, $\bfzeta_x^T = Z^T \X$.\\
(If the data matrices are scaled as in \cite{DDMD-SISC-2018}, then the scaling must be accordingly applied to $\bfzeta_x^T$.)
\end{proposition}
For the readers convenience, the case of canonical basis functions (the DMD framework) is summarized in Algorithm \ref{ALG:KS-DMD-1}. (Note that the algorithm is written in an EDMD framework, to express the DMD as a special case of the EDMD. A usual DMD-style formulation works on $\X$ and $\Y$, and the Rayleigh quotient is $\widehat{U}^T$.)
\begin{algorithm}
\caption{\label{ALG:KS-DMD-1} $(T, Z, \bfzeta_x^T) = \textsf{Koopman\_Schur\_SSMD}(\X, \Y)$}
	\begin{algorithmic}[1]
		%
		\STATE Compute the SVD $\X^T = W\Sigma V^*$
		\COMMENT{This step includes numerical rank determination.}
		\STATE $\widehat{U} = \Sigma^{-1}W^* (\Y^T V)$ \COMMENT{The Rayleigh quotient. See Section \ref{SS=DMD-connection-1}}
		\STATE Compute the Schur decomposition $\widehat{U} = Q T Q^*$ ; 
		\STATE $Z = VQ$ ; \COMMENT{The Schur vectors.}
		\STATE $\bfzeta_x^T = Q^T\Sigma W^T$ \COMMENT{The Schur functions values tabulated at the snapshots (\ref{zetax-SVD}).}
	\end{algorithmic}
\end{algorithm}

\subsubsection{Snapshot representation using a subset of modes}\label{SSS:SnapReprSubsetModes}
High fidelity representation of the data snapshots using only a selected subset of the computed modes allows for a model order reduction and revealing of underlying structures in the dynamics (e.g. identification and tracking of coherent structures in a fluid flow). The sparsity promoting DMD \cite{Jovanovic-Schmid-Nichols-2012}, \cite{jovschnicPOF14} is such a technique.

Our goal is to develop a similar snapshot reconstruction scheme, but for the Koopman-Schur framework. 
We start with a general result on a structured least squares problem that  generalizes a Khatri-Rao product based formulation  of the snapshot reconstruction problem  \cite{LS-Vand-Khatri-Rao-2020}.

\begin{theorem}\label{TM:LS-alpha}
Let $\X\in\Cz^{n\times m}$, $B\in\Cz^{n\times\ell}$, $C\in\Cz^{m\times\ell}$, 
$\Pond\in\Cz^{m\times m}$.
If the Khatri-Rao product $C\odot B$ is of full column rank and $\Pond$ is nonsingular, then the optimal solution $\bfalpha=(\alpha_1,\ldots, \alpha_\ell)$ of the least squares problem 
\begin{equation}\label{eq:LS-alpha}
    \| ( \X - B \left(\begin{smallmatrix} \alpha_1 &  & \cr
    & \cdot &  \cr
    &  & \alpha_\ell\end{smallmatrix}\right) C^T ) \Pond\|_F \longrightarrow \min_{\alpha_1,\ldots,\alpha_\ell}
\end{equation}
is unique and it can be compactly written as \begin{equation}\label{eq:alpha-solution} \bfalpha =  [(B^* B) \circ {(C^* \Pond^{T*} \Pond^T C)}]^{-1} [{C^*}\circ (B^* \X \Pond \Pond^*)]\ones,\;\; \ones = (1, \ldots, 1)^T\in\R^m ,
\end{equation}
where $\circ$ denotes the Hadamard matrix product.
\end{theorem}
{\sc Proof:}
The  linear least squares problem  (\ref{eq:LS-alpha}) can be rewritten into  
$$
vec(( \X - B \left(\begin{smallmatrix} \alpha_1 &  & \cr
    & \cdot &  \cr
    &  & \alpha_\ell\end{smallmatrix}\right) C^T ) \Pond) = 
    ( \Pond^T \otimes \Id_n) ( vec(\X) - (C\odot B)\left(\begin{smallmatrix} \alpha_1  \cr
\vdots \cr
\alpha_\ell\end{smallmatrix}\right) ),
$$
where $\otimes$ and $\odot$ denote the Kronecker and the Khatri-Rao matrix products, respectively, and $vec(\cdot)$ is the column-wise reshaping isomorphism between the vector spaces $\Cz^{n\times m}$ and $\Cz^{n\cdot m}$.
Further, if we define  
$$
S_\pond = (\Pond^T\otimes\Id_n)(C\odot B) = (\Pond^T C)\odot B ,
$$
then the minimization problem (\ref{eq:LS-alpha}) reads
$$
\| S_w \left(\begin{smallmatrix} \alpha_1  \cr
\vdots \cr
\alpha_\ell\end{smallmatrix}\right) - (\Pond^T\otimes\Id_n)vec(\X)\|_F\longrightarrow\min_{\alpha_1,\ldots,\alpha_\ell} .
$$
The assumptions of the theorem guarantee that $S_\pond$ has full column rank and the unique solution can be expressed using the Moore-Penrose pseudo-inverse as $\bfalpha = S_\pond^\dagger (\Pond^T\otimes\Id_n) vec(\X)=(S_\pond^* S_\pond)^{-1}S_\pond^*(\Pond^T\otimes\Id_n) vec(\X)$. The factors in this formula can be expressed, using the well known relations between the Hadamard, the Kronecker and the Khatri-Rao matrix products, as 
\begin{eqnarray}
S_\pond^* S_\pond &=& ((\Pond^T C)\odot B)^* ((\Pond^T C)\odot B) = [(\Pond^T C)^*(\Pond^T C)]\circ (B^* B) \nonumber \\ 
&=& {(C^*\Pond^{T*} \Pond^T C)}\circ (B^* B)  = (B^*B) \circ {(C^*\Pond^{T*}\Pond^T C)}\\
S_\pond^* (\Pond^T\otimes\Id_n)vec(\X) &=& (C\odot B)^*(\Pond^{T*}\Pond^T\otimes\Id_n)vec(\X) = [{C^*}\circ (B^* X \Pond \Pond^*)]\ones \nonumber 
\end{eqnarray}
which yields (\ref{eq:alpha-solution}).
\hfill$\boxtimes$
\begin{remark}
{\em
The condition that $C\odot B$ is of full column rank is satisfied if e.g. one of the two matrices is of full column rank and the other one has no zero column. This is a consequence of a Schur product theorem, see e.g. \cite[Section 5.2]{hor-joh-91}.
}
\hfill $\boxtimes$\end{remark}
\begin{remark}
{\em
Clearly, with an additional nonsigular matrix $\Pond_1$, the problem can be formulated as 
\begin{equation}\label{eq:LS-alpha-2}
    \|\Pond_1 ( \X - B \left(\begin{smallmatrix} \alpha_1 &  & \cr
    & \cdot &  \cr
    &  & \alpha_\ell\end{smallmatrix}\right) C^T ) \Pond\|_F \longrightarrow \min_{\alpha_1,\ldots,\alpha_\ell}
\end{equation}
with the solution
$\bfalpha =  [(B^*\Pond_1^*\Pond_1 B) \circ {(C^* \Pond^{T*} \Pond^T C)}]^{-1} [{C^*}\circ (B^* \Pond_1^*\Pond_1 \X \Pond \Pond^*)]\ones.$
In this case $\Pond_1$ can be e.g. the Cholesky factor of the inverse noise covariance matrix. Throughout this paper, for the sake of brevity, we have $\Pond_1=\Id_n$.
See \cite{LS-Vand-Khatri-Rao-2020} for more details.
}
\hfill $\boxtimes$\end{remark}
The formula (\ref{eq:alpha-solution}) is in  particular interesting in our case where
$C=\bfzeta_x$ and $B=\Xi^T=Z^{*T}$ is orthonormal and thus $B^*B=\Id_\ell$. 
\begin{corollary}\label{COR:alpha:B-orth} If in Theorem \ref{TM:LS-alpha} the matrix $B$ is orthonormal ($B^*B=\Id_\ell$), then the
solution (\ref{eq:alpha-solution}) of the least squares problem (\ref{eq:LS-alpha}) reads
$$
\bfalpha= \mathrm{diag}(1/\|\Pond^T C(:,j)\|_2^2)_{j=1}^\ell [{C^*}\circ (B^* \X \Pond\Pond^* )]\ones.
$$
\end{corollary}
\begin{remark}
{\em 
The coefficients in Corollary \ref{COR:alpha:B-orth} are actually the Fourier coefficients and  the formula for $\bfalpha$ follows directly from the projection theorem. Namely, if $b_i$ and $c_i$ denote the $i$th column of $B$ and $C$, respectively, then the approximation problem reads
$$
\X\Pond \approx \sum_{i=1}^\ell \alpha_i b_i c_i^T\Pond,
\;\;\mbox{where}\;\;\langle b_i c_i^T\Pond, b_jc_j^T\Pond \rangle_F = \bfdelta_{ij} \|b_i\|_2^2 \|\Pond^T c_i\|_2^2.
$$
Here $\bfdelta_{ij}$ is the Kroneecker symbol and $\langle \cdot,\cdot\rangle_F$ is the Frobenius inner product. The optimal coefficients are thus
$$
\alpha_i\! =\! \frac{\langle \X\Pond,b_ic_i^T\Pond\rangle_F}
{\langle b_ic_i^T\Pond, b_ic_i^T\Pond\rangle_F} \!=\! \frac{\mathrm{Trace}(\Pond^*c_i^{T*}b_i^* \X\Pond)}{ \|b_i\|_2^2 \|\Pond^T c_i\|_2^2} \!=\!
 \frac{\mathrm{Trace}(b_i^* \X\Pond\Pond^*c_i^{T*} )}{ \|\Pond^T c_i\|_2^2} \!=\!
  \frac{[c_i^* \!\circ \! (b_i^* \X\Pond\Pond^*)]\ones}{ \|\Pond^T c_i\|_2^2},
$$
precisely as stated in Corollary \ref{COR:alpha:B-orth}.
}
\hfill $\boxtimes$\end{remark}
The role of the matrix $\Pond$ can be to weigh the importance  of individual snapshots in the overall reconstruction error. This is achieved by choosing $\Pond$ as diagonal positive definite.\footnote{Actually, $\Pond$ can be semidefinite, where $\pond_i=0$ means that the $i$th snapshot is simply excluded. After removing all zero-weighted snapshots, we proceed with the problem with positive weights.}
\begin{corollary}
 Let $B=\Xi^T(:,[i_1,\ldots, i_{\ell}])=(\eta_{i_1},\ldots, \eta_{i_\ell})$,  $C=\bfzeta_x(:,[i_1,\ldots, i_{\ell}])$, and $\Pond=\mathrm{diag}(\pond_i)_{i=1}^m$, $\pond_i>0$ for all $i=1,\ldots, m$. Then the optimal reconstruction $\X\approx B\bfalpha C^T$ in the sense of Theorem \ref{TM:LS-alpha} is obtained with
$$
\bfalpha = \mathrm{diag}(\frac{1}{\sum_{k=1}^m \pond_k^2|\zeta_{i_j}(\x_k)|^2})_{j=1}^\ell \;
[({\underbrace{(\overline{\zeta_{i_j}(\x_k)})_{j,k=1,1}^{\ell,m}}_{C^*}})\circ (Q(:,[i_1,\ldots,i_\ell])^T V^T \X \Pond^2)]\ones .
$$
The computational complexity of the above formula (given its ingredients) is independent of the state space dimension $n$. In a Koopman-Schur-ESSMD setting, the kernel trick is used, as outlined in Section \ref{SS:eval-Schur-f}.
If $\bfPsi_x=\X^T$, as in  Proposition \ref{PROP:snap-rep-2}, then
$\bfalpha = (1,\ldots,1)^T$.
\end{corollary}
{\sc Proof}: Note that 
$
B^*\X = (Z^{T*}(:,[i_1,\ldots,i_\ell]))^* X = Q(:,[i_1,\ldots,i_\ell])^T V^T \X .
$
If $\bfPsi_x=\X^T$, then we have (from the SVD of $\bfPsi_x$)  the SVD of $\X$ as well, and
$$
B^* \X = Q(:,[i_1,\ldots,i_\ell])^T V^T V^{*T}\Sigma W^T=Q(:,[i_1,\ldots,i_\ell])^T\Sigma W^T = \bfzeta_x(:,[i_1,\ldots,i_\ell])^T .
$$
Since $B^*B=\Id_\ell$, Corollary \ref{COR:alpha:B-orth} yields
$$
\bfalpha = \mathrm{diag}(\frac{1}{\sum_{k=1}^m \pond_k^2|\zeta_{i_j}(\x_k)|^2})_{j=1}^\ell \;
[({\underbrace{(\overline{\zeta_{i_j}(\x_k)})_{j,k=1,1}^{\ell,m}}_{C^*}})\circ (\underbrace{Q(:,[i_1,\ldots,i_\ell])^T\Sigma W^T}_{C^T}\Pond^2)]\ones .
$$
It is easily checked that the above expression actually reads $\bfalpha = (1,\ldots,1)^T$. \hfill $\boxtimes$

%

\subsubsection{Reordering trick}

A reduced order representation of the snapshots may be based on a selection of eigenvalues $\lambda_{i_1} = T_{i_1i_1}, \ldots, \lambda_{i_{\ell}}=T_{i_\ell i_\ell}$. In the standard DMD/EDMD approaches, one simply takes the corresponding modes and solves a least squares problem for the coefficients. In the Koopman-Schur setting, the basis for the representation is an orthonormal basis that spans the spectral subspace encompassing the selected eigenvalues. In a Schur form, such a basis is in the leading columns of the unitary factor, provided the eigenvalues are at the corresponding leading positions (upper left corner) of the triangular factor; the decomposition is then truncated and only selected values are used in the next steps. Another situation where reordering is needed is when spurious eigenvalues appear (e.g. that are far outside unit circle) -- reordering can move them to the bottom right part of the triangular form and then truncated from the decomposition.

The reordering trick, that enables orthonormal representation with any subset of the eigenvalues (Section \ref{SSS:SnapReprSubsetModes}) is to reorder the Schur form (Section \ref{SS=ordered-Schur}) so that in the new triangular Schur factor $\widetilde{T}$ we have $\widetilde{T}_{11}=\lambda_{i_1}, \ldots, \widetilde{T}_{\ell\ell}=\lambda_{i_\ell}$, and then to truncate it at index $\ell$. 
%
\subsection{Using the Schur functions for forecasting}\label{SS=KS_forecasting}
The Koopman mode decomposition is a powerful tool for revealing coherent/latent structure of the dynamics under study, model order reduction and for forecasting the future values of the observables of interest. For the latter, the key is the accuracy of the numerically computed eigenfunctions and the eigenvalues. 

The proposed Koopman-Schur decomposition can also be adapted for forecasting. 
The forecasting mechanism works as follows.
Once we have the representation (\ref{eq:g-repr}), the relation (\ref{eq:Uzeta=zetaT(x)}) implies 
$$
\begin{pmatrix} (\KO g_1)(x) & \ldots & (\KO g_\ell)(x)\end{pmatrix} \backsimeq \begin{pmatrix}\zeta_1(x) & \ldots & \zeta_r(x) \end{pmatrix} T \Xi , 
$$
and, inductively, 
\begin{equation}\label{eq:Ukg}
\begin{pmatrix} (\KO^k g_1)(x) & \ldots & (\KO^k g_\ell)(x)\end{pmatrix} \backsimeq \begin{pmatrix}\zeta_1(x) & \ldots & \zeta_r(x) \end{pmatrix} T^k \Xi,\;\; k=1, 2, \ldots. 
\end{equation}
Hence, pushing the observables $g_i$ forward ($\KO^k g_i = g_i\circ \DDS \circ \cdots \circ\DDS$) is, in the representation (\ref{eq:g-repr}), realized by raising the power of the triangular factor $T$. This includes the forecasting, i.e. extrapolating the representation of the present snapshot into the future.
In terms of the values at the snapshots (see (\ref{eq:GT-reconstruct})) the relations (\ref{eq:Ukg}) read 
\begin{equation}\label{eq:Ukg-at_xj}
((\KO^k g_i)(\x_j))_{i,j=1}^{\ell,m} \approx \Xi^T (T^T)^k\bfzeta_x^T,\;\; k=1, 2, \ldots.
\end{equation}
Of course, with repeated applications of $\KO$, the residual will build up (see (\ref{eq:residual-rhoTZ})) and the approximations (\ref{eq:Ukg}) will deteriorate in the sense of $\backsimeq$ .
Nevertheless, we can use (\ref{eq:Ukg-at_xj}) to approximate, for $j=2, 3, \ldots$
$$
\Xi^T \!\begin{pmatrix}
\zeta_1(\x_{j})\cr \zeta_2(\x_{j})\cr\vdots\cr \zeta_r(\x_{j}) 
\end{pmatrix} \!\!\approx\!\!
\begin{pmatrix}
g_1(\x_j)\cr g_2(\x_j)\cr\vdots\cr g_\ell(\x_j) 
\end{pmatrix} \!\!=\!\!
\begin{pmatrix}
(g_1\circ \DDS)(\x_{j-1})\cr (g_2\circ\DDS)(\x_{j-1})\cr\vdots\cr (g_\ell\circ\DDS)(\x_{j-1}) 
\end{pmatrix} \!\!=\!\!
\begin{pmatrix}
(\KO g_1)(\x_{j-1})\cr (\KO g_2)(\x_{j-1})\cr\vdots\cr (\KO g_\ell)(\x_{j-1}) 
\end{pmatrix}
\!\approx\! \Xi^T T^T\!\begin{pmatrix}
\zeta_1(\x_{j-1})\cr \zeta_2(\x_{j-1})\cr\vdots\cr \zeta_r(\x_{j-1}) 
\end{pmatrix} ,
$$
i.e. 
\begin{equation}\label{eq:consistent-step}
   \begin{pmatrix}
\zeta_1(\x_{j})\cr \zeta_2(\x_{j})\cr\vdots\cr \zeta_r(\x_{j}) 
\end{pmatrix} \approx  T^T\!\begin{pmatrix}
\zeta_1(\x_{j-1})\cr \zeta_2(\x_{j-1})\cr\vdots\cr \zeta_r(\x_{j-1}) 
\end{pmatrix}, \;\;j=2, 3, \ldots, m, m+1.
\end{equation}
The quality of the approximation $\approx$ in (\ref{eq:consistent-step}) is an indirect indicator of the success/failure of the entire discretization process from Section \ref{SS=Num-Compress}.

Recall that for an arbitrary $x\in\Cz^n$
$$
\begin{pmatrix}
\zeta_1(x)\cr \zeta_2(x)\cr\vdots\cr \zeta_r(x) 
\end{pmatrix} = Q^T \Sigma^{-1}W^T
\begin{pmatrix}
\kerf(x,\x_1)\cr \kerf(x,\x_2)\cr \vdots\cr \kerf(x,\x_m)
\end{pmatrix} = Q^T\Sigma^{-1}W^T \X^* x = \underbrace{Q^T\Sigma^{-1}W^T W^{*T}\Sigma V^T}_{Z^T} x ,
$$
where the last two equalities hold in the case $\bfPsi_x=\X^T$, $\kerf(x,y)=\langle x, y\rangle_{\Cz^n}=y^* x$. In particular, for  $x=\x_{m+1}$, 
$$
\begin{pmatrix}
\zeta_1(\x_{m+1})\cr \zeta_2(\x_{m+1})\cr\vdots\cr \zeta_r(\x_{m+1}) 
\end{pmatrix}\! = Q^T \Sigma^{-1}W^T\!
\begin{pmatrix}
\langle \x_{m+1},\x_1\rangle_{\Cz^n}\cr \langle \x_{m+1},\x_2\rangle_{\Cz^n}\cr \vdots\cr \langle \x_{m+1},\x_m\rangle_{\Cz^n}
\end{pmatrix} 
= \overline{Z}^* \x_{m+1}=Z^T \x_{m+1},
$$
as it should be from the least squares approximation theory (Proposition \ref{PROP:snap-rep-2}).
Note that (\ref{eq:consistent-step}) represents an explicit time stepping that can be deployed for the purposes of forecasting beyond the index $j=m+1$.

\subsection{The real Koopman-Schur SSMD}\label{SS=Real-KS-SSMD}
If the data is real, it is advantageous to perform as much as possible of the computation using only the real computer arithmetic.  Since the eigenvalues are in general complex, 
we can reformulate the algorithm using the real Schur form \cite[Sec. 2.3]{hor-joh-90}. This means that in the Schur form (\ref{eq:Schur-hat-U}) of a real $\widehat{U}$ the matrix $Q$ is real orthogonal and $T$ is block upper triangular with diagonal blocks $1\times 1$ for each real eigenvalue and  $2\times 2$ for each complex conjugate pairs of eigenvalues. For instance, this real form is the default in the Matlab function \texttt{schur} and the complex form is obtained with the input option \texttt{'complex'}. The LAPACK subroutines \texttt{SGEES, DGEES} are designed to compute the real Schur form of real matrices. 

Using the real Schur form for real data is recommended and it is implemented in our software. The required changes (as compared with the complex form) are merely technical and require most additional effort during software development.
Some care is needed when defining the nested invariant subspaces (see (\ref{eq:U-Schur-i})) in the sense that complex conjugate pairs should always be taken together, with the corresponding Schur vectors.
Further, if a reordering of the Schur form is needed, the reordering algorithm should be adapted to the real form.
Our LAPACK-based implementation of the Koopman-Schur SSMD method (a separate work) is designed to use the real Schur form for real data, similarly as we have recently done for DMD implementations \cite{lawn-298-drmac}, \cite{lawn-300-drmac}.

\subsection{QR compressed Koopman-Schur SSMD}
The QR compressed DMD \cite{DDMD-SISC-2018}, with implementation details \cite[Section 3.7]{lawn-298-drmac} can considerably increase the run time efficiency of the DMD in case of high dimensional data snapshots, and it can be used for fast updates (adding new data snapshots, single or in batches) and down-dates (discarding oldest data snapshots, while receiving new ones) for the streaming DMD. The technique is easily transferred to the Hermitian DMD
in \cite{lawn-300-drmac}. 

It is possible to apply the QR compression to the KS-SSMD Algorithm \ref{ALG:KS-DMD-1} and to formulate a streaming version. We will provide details in a forthcoming work.

\section{Numerical evaluation}\label{S=NumericalEvaluation}
We now illustrate the Koopman-Schur framework using numerical examples. Software prototypes of the presented algorithms are developed using Matlab and the numerical examples are selected to check the validity and the functionality of the new proposed approach. 
The experiments provide insights and show differences between the methods, both row-wise (DMD vs. EDMD) and column-wise (diagonalization vs. Schur form) in Table \ref{Table-1}.
\begin{table}[H]
\begin{tabular}{|l|l|l||}
\hline\hline
 &  Dynamic Mode Decomp. & Extended Dynamic Mode Decomp.\\ \hline
 diagonalization & DMD & EDMD \\
 Schur form  & Koopman-Schur-SSMD & Koopman-Schur-ESSDMD\\ \hline\hline
\end{tabular}
\caption{The four methods used in the numerical experiments in this section. In the EDMD-based methods (second column) the kernel-trick formulation is used.\label{Table-1}}
\end{table}

We use two kernel functions
\begin{enumerate}
\item $k_1(x,y) = y^* x = \langle x, y\rangle_{\Cz^n}$. With this choice EDMD and DMD are theoretically  equivalent 
but EDMD is numerically more difficult because of higher condition numbers of the matrices used in the algorithm.
\item $k_2(x,y) = e^{-\frac{\|x-y\|_2^2}{2\sigma^2}}$. The Gaussian kernel is popular and powerful kernel function, and for good performance proper choice of the hyper-parameter $\sigma$ is critical. We will not tackle the question of choosing $\sigma$; instead we have empirically selected a value that works well on our examples.
\end{enumerate}

\subsection{What is being tested}\label{SSS=WhatisTested}
\noindent The main goal of the experiments presented here is to illustrate the proposed framework. The themes of the test are as follows:

\paragraph{Eigenvalues.} The eigenvalues computed by all four methods are compared. 
If in some test case we have the eigenvalues of the matrix that generates the data, we use them
as ground truth reference values. Such synthetic examples are useful for the first phase of implementation and
debugging of the code. For the data from simulations of physical processes we study the numerical behaviour of the computed eigenvalues along a sequence of sliding data windows.

\paragraph{Snapshot reconstruction.} The supplied data snapshots are represented using the Schur functions (see Section \ref{SS=Observ_Representation_Schur_f}) and the residuals of the representation are computed.
At each time step, the maximal reconstruction error over all data snapshots in the active windows is computed and recorded. At the end of each test, the maximal errors over all windows are displayed for all four methods.  

\paragraph{Consistency.} We use relation (\ref{eq:consistent-step}) as an implicit test of
          the correctness of the implementation by computing the residuals 
$$
\left\| \left(\begin{smallmatrix}
\zeta_1(\x_{j}) \cr \zeta_2(\x_{j}) \cr \vdots \cr  \zeta_r(\x_{j}) 
\end{smallmatrix}\right) - T^T\! \left(\begin{smallmatrix}
\zeta_1(\x_{j-1}) \cr \zeta_2(\x_{j-1}) \cr \vdots \cr \zeta_r(\x_{j-1}) 
\end{smallmatrix}\right)\right\|_2, \;\;j=2, 3, \ldots, m, m+1. 
$$
These residuals check whether, with given data, the data driven compression of the Koopman operator, outlined in  Section \ref{SS=Num-Compress}, was numerically computed to sufficient accuracy. Further, if these residuals are large (meaning that the known past and the present cannot be reproduced) then the forecasting is doomed to fail.
Maximal residuals are computed for each active window and displayed for the methods in the second row of Table \ref{Table-1}.

\paragraph{Forecasting skill.} Selected data sets are used to show that the new method has the prediction skill as the standard DMD and EDMD implementations. \\

\noindent We should also emphasize what is not included in this presentation of the new method. High performance software implementations, numerical robustness in terms of the core numerical linear algebra,  and adaptations to particular applications (such as dynamic adaptation for forecasting dynamics with rapid changes, or revealing coherent structures) are challenging themes that  are the subjects of separate work. A LAPACK style implementation that addresses more numerical details, following our recent work \cite{lawn-298-drmac}, \cite{lawn-300-drmac}, is in progress and will be soon available.

\subsection{Examples}
We have tested the new proposed scheme on many examples. Here we present the result of two test cases and discuss the
key elements of the experiment outlined in
Section \ref{SSS=WhatisTested}.

\begin{example}\label{EX:cylinder}
{\em
 In the first example, we use a simulated data of a 2D laminar flow around a cylinder. 
 The $1201$ data snapshots of dimension $201\times 101$ are reshaped into a $20301\times 1201$ matrix and used as test data, divided into $\X$ and $\Y$ as explained in Section \ref{S=Intro}.
The sliding data window is set of fixed size $w=100$, and simulation is run for $100$ time steps. At each step a new snapshot is received, the oldest snapshot is discarded and the DMD, EDMD, Koopman-Schur-SSMD and Koopman-Schur-ESSMD decompositions are computed.  The Koopman-Schur-SSMD and the Koopman-Schur-ESSMD are used for forecasting, to check the framework outlined in Section \ref{SS=KS_forecasting}. The forecasting  scheme is kept simple, following the formulas derived in Section \ref{SS=KS_forecasting}, and no attempt is made to dynamically resize the active window. The forecasting horizon is set to $\tau=40$.
Both kernel functions are tested, and for the Gaussian kernel the parameter $\sigma$ is set to $300$ without any fine tuning. 

Figure \ref{FIG:cyl:k1w100-recon-cons} and Figure \ref{FIG:cyl:k2w100-recon-cons} show the reconstruction errors for all four methods and the consistency errors for the Koopman-Schur based methods, respectively. 
The eigenvalues computed by all four methods and the condition number of the key matrices, using the kernel $k_1(\cdot,\cdot)$, are shown in Figure \ref{FIG:cyl:k1w100-eig-cond}. The prediction skills of the KS-SSMD and KS-EDMD (with the kernels $k_1(\cdot,\cdot)$ and $k_2(\cdot,\cdot)$ ) are shown in Figure \ref{FIG:cyl-k1-k2-pred}. The forecast values are computed by direct implementation of the formulas from Section \ref{SS=KS_forecasting}, without any additional modifications (such as dynamic window resizing, removing spurious modes, curbing ill-conditioning etc.), and the results are encouraging and promising.

Figure \ref{FIG:cylk2w100eigscond} shows the eigenvalues of all four methods as computed from the last active window (with the kernel $k_2(\cdot,\cdot)$) and the condition numbers of the relevant matrices, As expected, there are two groups of the eigenvalues: the one computed by DMD and KS-SSMD and the other by EDM and KS-ESSMD. Since the two groups of methods operate using different function dictionaries, they not necessarily return the same eigenvalues. But if the values from those two groups get close together, it is an indication of numerical convergence. We do not have analytical result but empirical evidence is interesting, We have collected all computed Ritz values (from all $100$ windows) ans plotted them all as in the left panel in Figure \ref{FIG:cylk2w100eigsall}. Then we identified in the left panel in Figure \ref{FIG:cylk2w100eigscond} the eigenvalues that had been computed nearly the same by all methods, and zoomed around them in the left panel of Figure \ref{FIG:cylk2w100eigsall}. Two particular selection are shown in the middle and the right panel. Tightly clustered values computed by different methods are easily identified by visual inspection. A more systematic empirical approach is to use, say, machine learning techniques to cluster and follow the eigenvalues computed from a sequence of sliding data windows,\footnote{We tackle this in a separate work.}
\begin{figure}[h]	
\includegraphics[width=3.in,height=2.in]{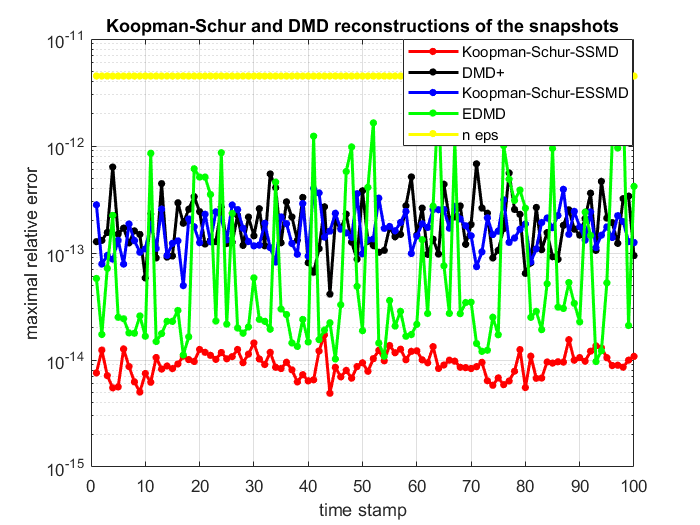}\hspace{-4mm}
\includegraphics[width=3.in,height=2.in]{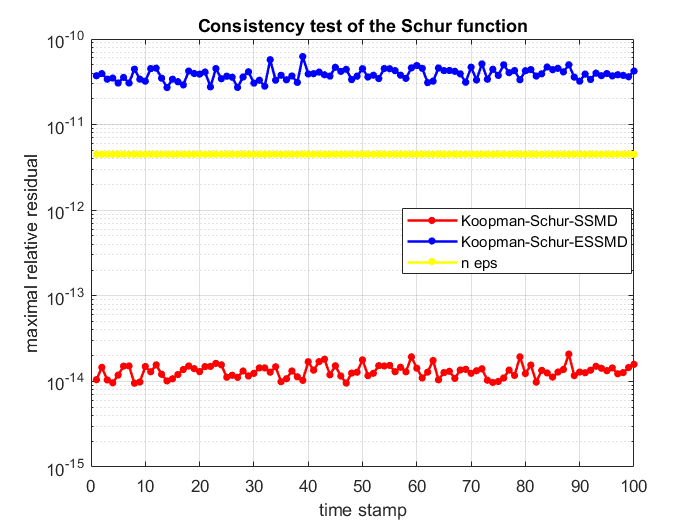}
\caption{(Example \ref{EX:cylinder}) The reconstruction error and the consistency of the computed Schur functions. The kernel function is $k_1(\cdot,\cdot)$.  The yellow horizontal line indicates the level of machine precision {\tt eps} times the state space dimension $n$. \label{FIG:cyl:k1w100-recon-cons}}
\end{figure}
\begin{figure}[h]
\includegraphics[width=3.in,height=2.in]{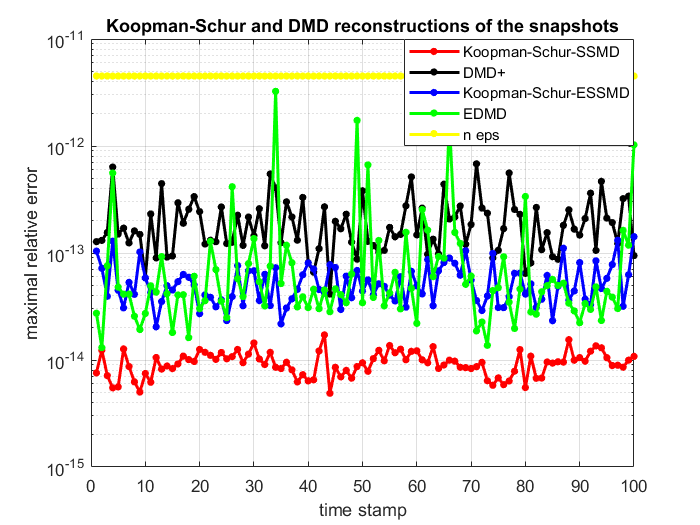}\hspace{-4mm}
\includegraphics[width=3.in,height=2.in]{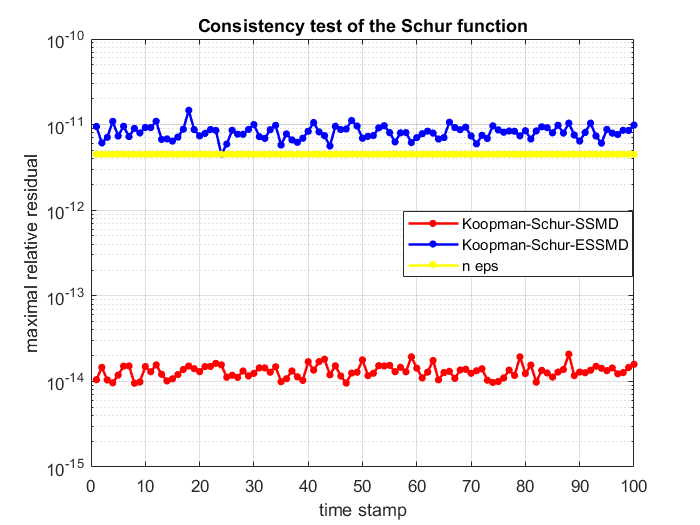}
\caption{(Example \ref{EX:cylinder}) The reconstruction error and the consistency of the computed Schur functions.  The kernel function is $k_2(\cdot,\cdot)$. The yellow horizontal line indicates the level of machine precision {\tt eps} times the state space dimension $n$.\label{FIG:cyl:k2w100-recon-cons}}
\end{figure}

\begin{figure}[h]	
\includegraphics[width=3.in,height=2.in]{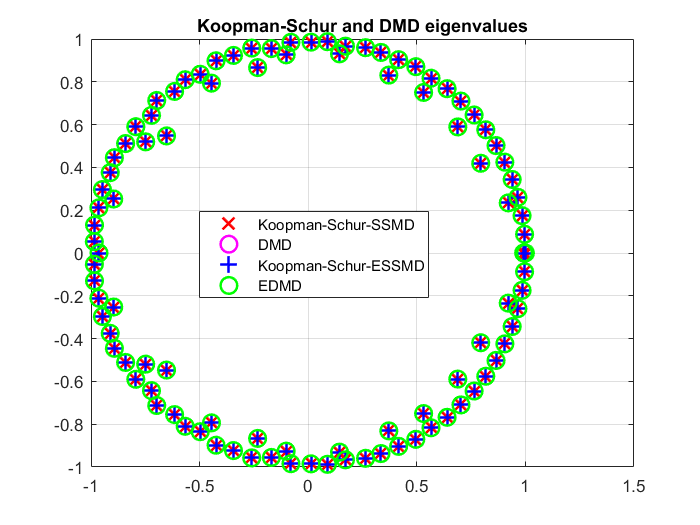}\hspace{-4mm}
\includegraphics[width=3.in,height=2.in]{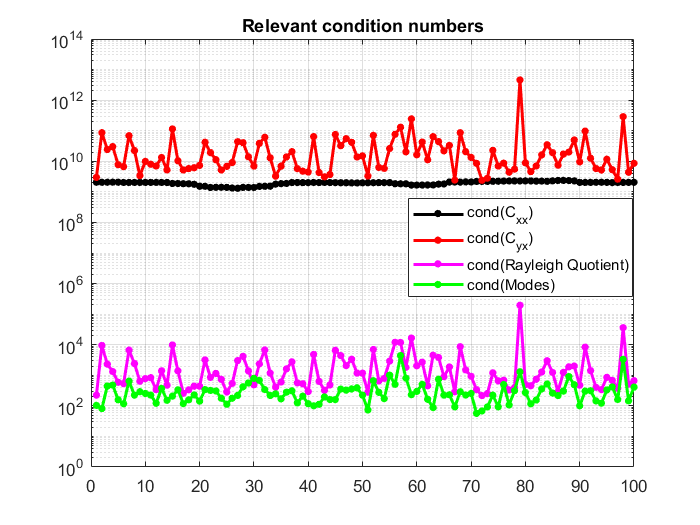}
\caption{(Example \ref{EX:cylinder}) \emph{Left}: The eigenvalues (Ritz values) computed in the last active window. \emph{Right}: The relevant condition numbers. The kernel function is $k_1(\cdot,\cdot)$. \label{FIG:cyl:k1w100-eig-cond}}
\end{figure}

\begin{figure}[h]	
\includegraphics[width=0.34\linewidth,height=1.4in]{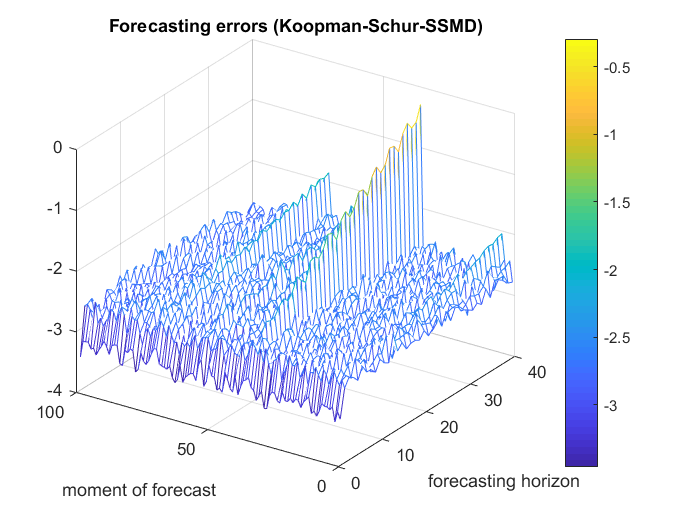}\hspace{-3mm}
\includegraphics[width=0.34\linewidth,height=1.4in]{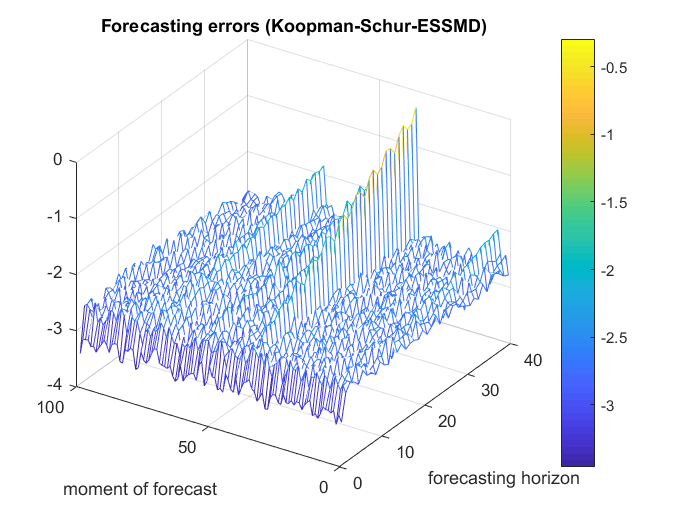}\hspace{-3mm}
\includegraphics[width=0.34\linewidth,height=1.4in]{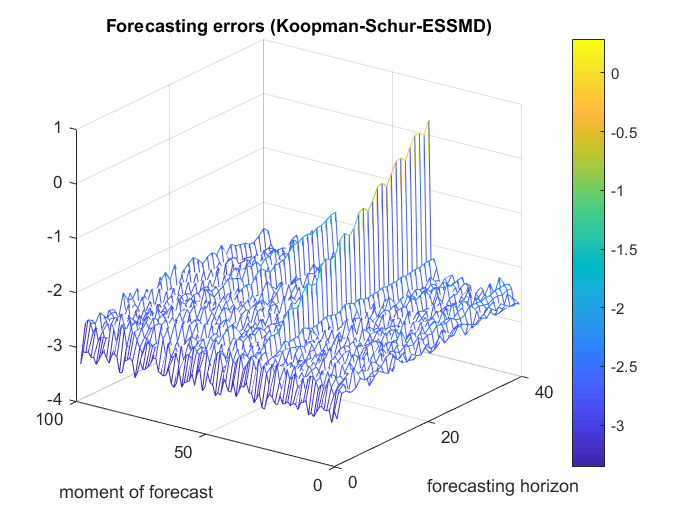}
\caption{(Example \ref{EX:cylinder}) The forecast error with the horizon $\tau=40$. \emph{Left}: Koopman-Schur-SSMD. \emph{Middle}: Koopman-Schur-ESSMD with the kernel function $k_1(\cdot,\cdot)$.
The KS-SSMD and the KS-ESSMD method are equivalent and the difference is only due to the finite precision arithmetic.
\emph{Right}: Koopman-Schur-ESSMD with the kernel function $k_2(\cdot,\cdot)$.\label{FIG:cyl-k1-k2-pred}}
\end{figure}

\begin{figure}[h]	
\includegraphics[width=3.1in,height=2.4in]{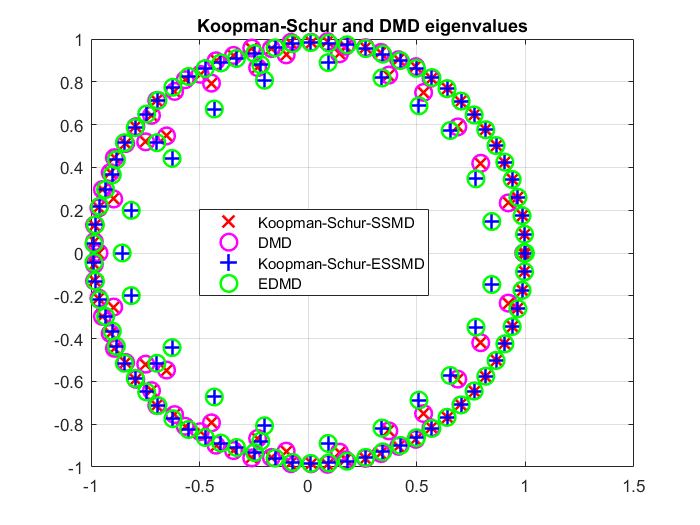}\hspace{-4mm}
\includegraphics[width=3.1in,height=2.4in]{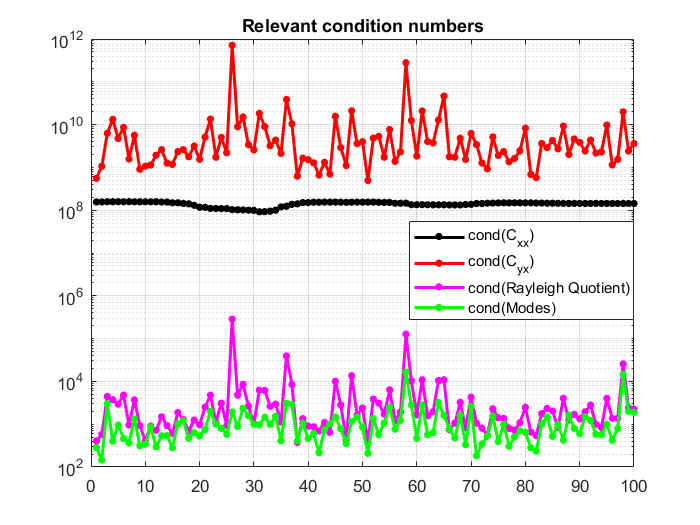}
\caption{(Example \ref{EX:cylinder}) \emph{Left}: The eigenvalues (Ritz values) computed in the last active window. \emph{Right}: The relevant condition numbers. The kernel function is $k_2(\cdot,\cdot)$. \label{FIG:cylk2w100eigscond}}
\end{figure}
\begin{figure}[h]	
\includegraphics[width=2.in,height=1.5in]{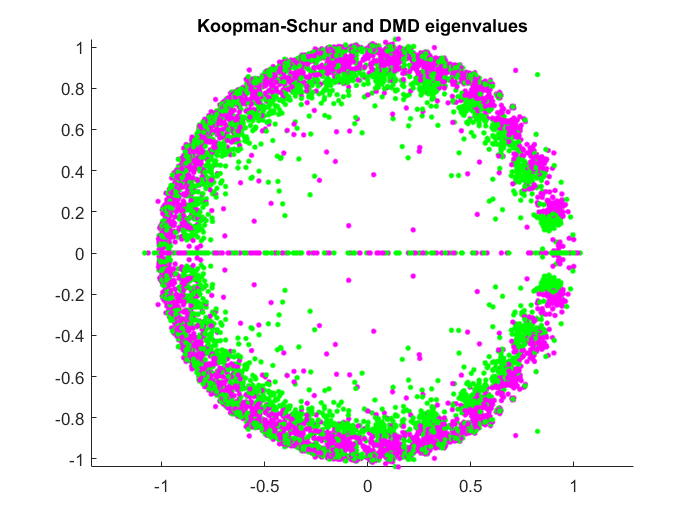}\hspace{-4mm}
\includegraphics[width=2.in,height=1.5in]{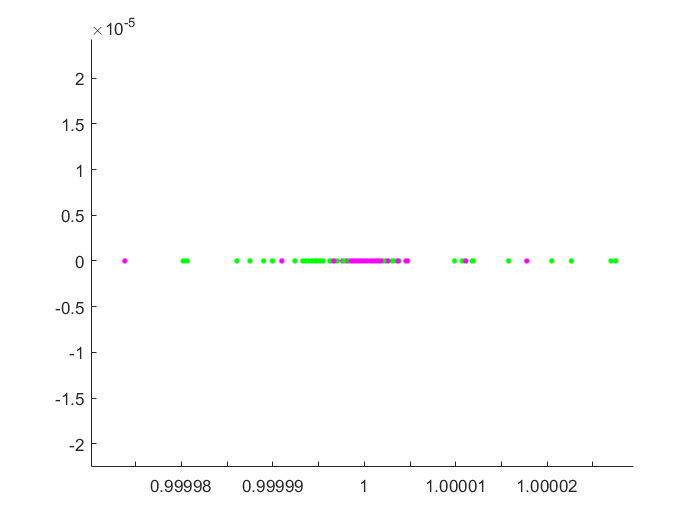}\hspace{-4mm}
\includegraphics[width=2.in,height=1.5in]{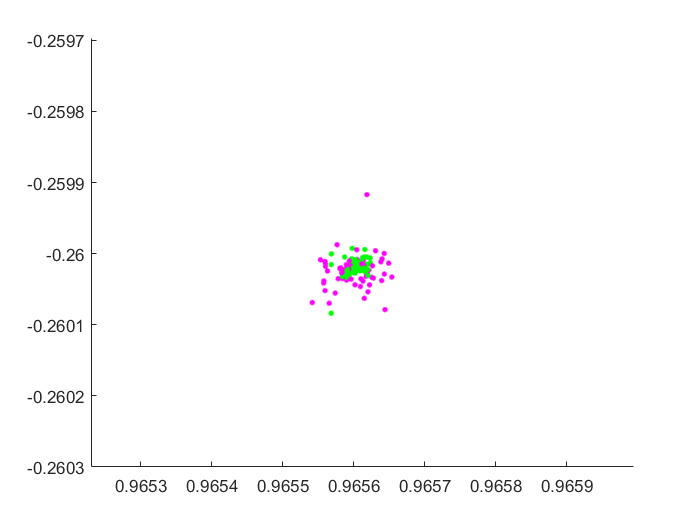}
\caption{(Example \ref{EX:cylinder}) \emph{Left}: All Ritz values computed throughout $100$ time steps.
\emph{Middle and right}: Zoomed in around values identified in Figure \ref{FIG:cylk2w100eigscond} (all four methods computing nearly identical values) to show how the sequence of Ritz values build tight clusters. This motivates using clustering method throughout iterations and techniques based e.g. on diffusion maps for tracking the Ritz values to reveal possible numerical convergence. (See Figure \ref{FIG:vort-eigs-70-100} below.) \label{FIG:cylk2w100eigsall}}
\end{figure}

\hfill$\boxtimes$
 }
\end{example}

\noindent Example \ref{EX:cylinder} is well-conditioned and it is a good test case to check the concept and the software implementation. It shows good prediction skill of the method based on the Schur decomposition.

\begin{example}\label{EX:vorticity}
{\em
The next example uses vorticity data of a quasi-two-dimensional Kolmogorov flow \cite{Tithof-etal-q2d-kolmog-2017}.
This is a difficult example, for discussions and details we refer to \cite{Suri-thesis-2017}, 
\cite{Tithof-etal-q2d-kolmog-2017}, \cite{PhysRevLett.118.114501}. (The methods tested here (Table \ref{Table-1})
are, as in Example \ref{EX:cylinder}, entirely data driven and completely oblivious to the underlying nature of the data.)
The experiment is performed with a similar setting as in Example \ref{EX:cylinder}.
Figure \ref{FIG:vort-recon-con-70-100} shows the reconstruction errors and the condition numbers of two runs, with the window sizes $w=70$ and $w=100$. The impact of the condition number of the matrix of the modes to the reconstruction of the snapshots is apparent.

 \begin{figure}[h]	 
 \includegraphics[width=3.1in,height=2.4in]{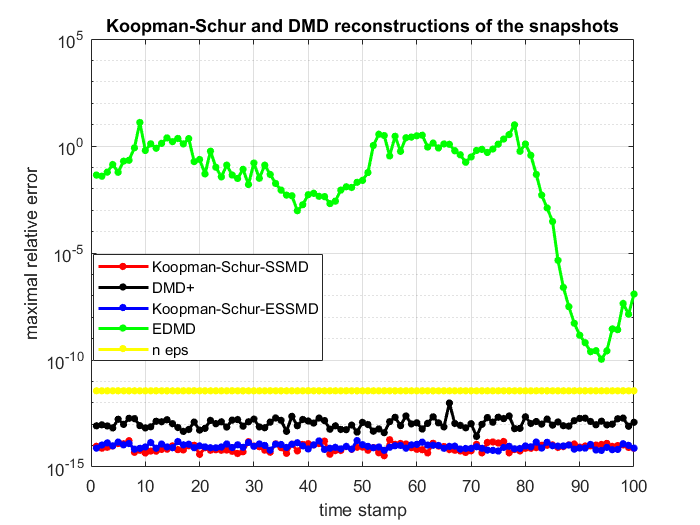}\hspace{-4mm}
\includegraphics[width=3.1in,height=2.4in]{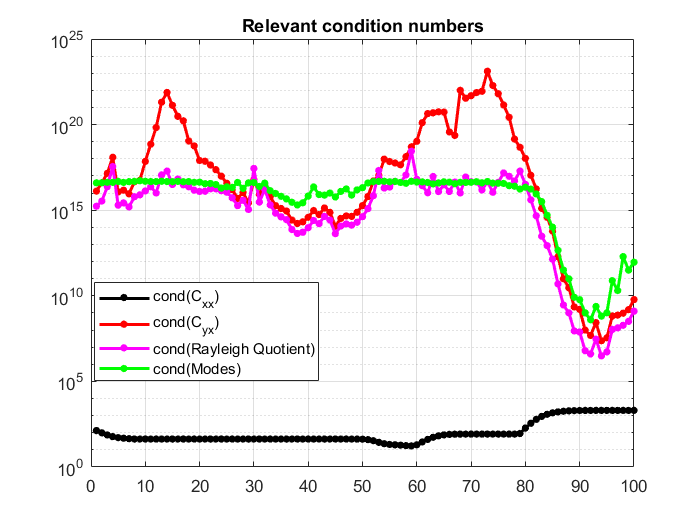}
\includegraphics[width=3.1in,height=2.4in]{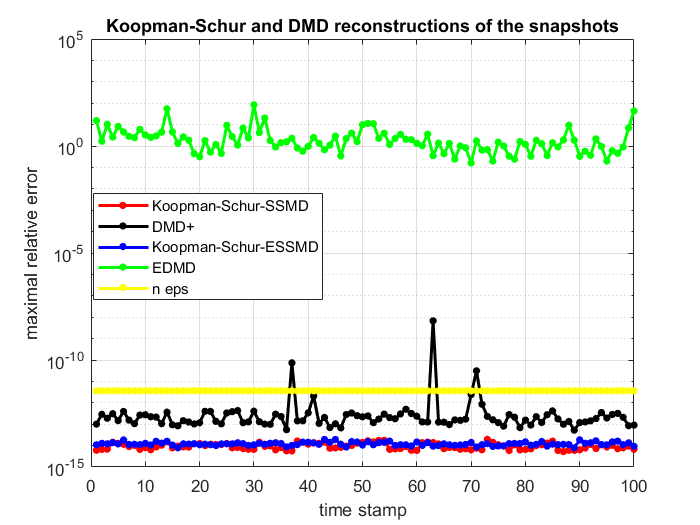}\hspace{-4mm}
\includegraphics[width=3.1in,height=2.4in]{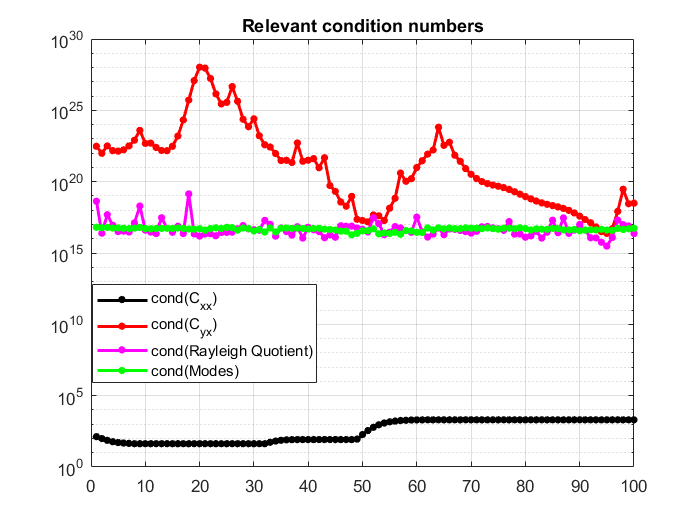}
\caption{(Example \ref{EX:vorticity}) \emph{First row}: Active window width $w=70$.
\emph{Second row}: Active window width $w=100$. The kernel function is $k_2(\cdot,\cdot)$. \label{FIG:vort-recon-con-70-100}}
\end{figure}

 \begin{figure}[h]	
\includegraphics[width=0.35\linewidth,height=1.5in]{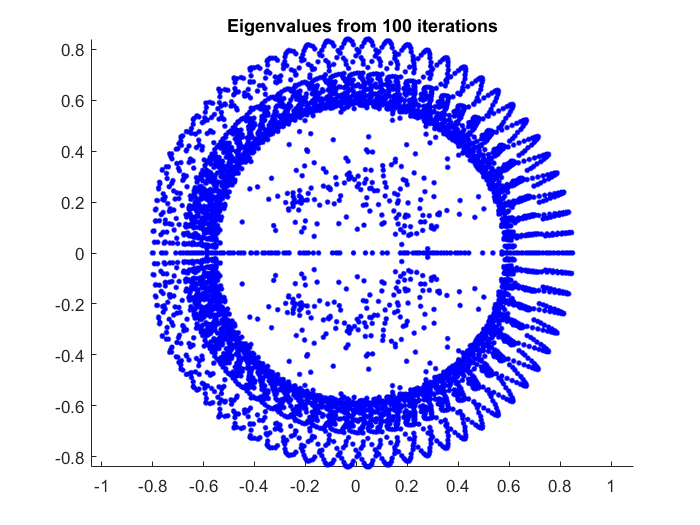}\hspace{-2mm}
\includegraphics[width=0.35\linewidth,height=1.5in]{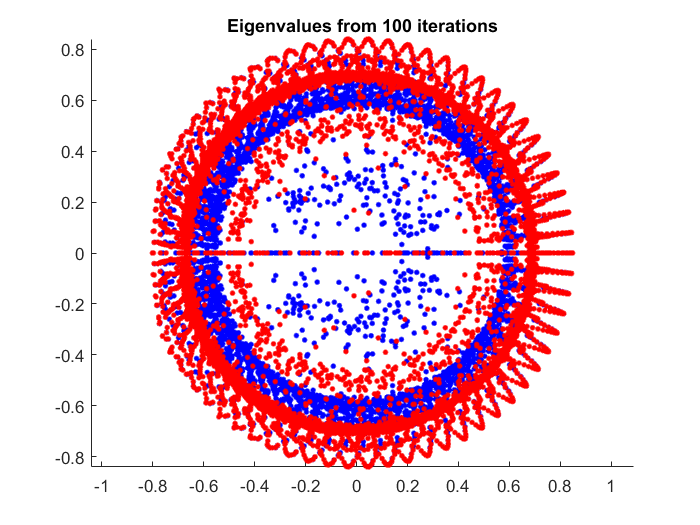}\hspace{-2mm}
\includegraphics[width=0.30\linewidth,height=1.5in]{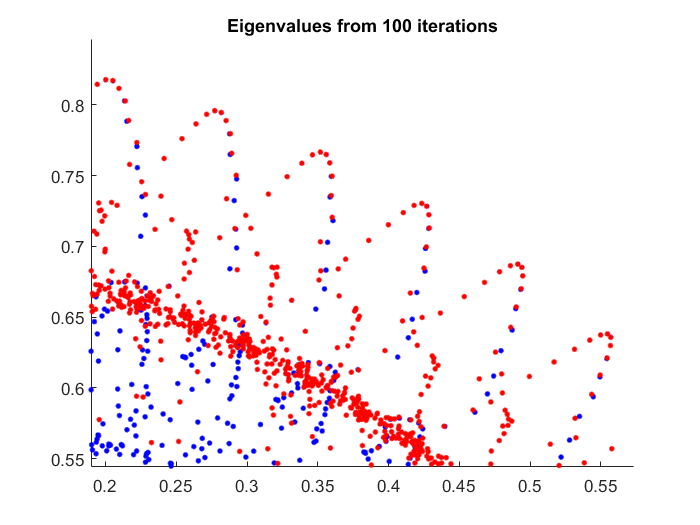}
\caption{(Example \ref{EX:vorticity}) The eigenvalues from $100$ sliding windows with the kernel $k_2$.
\emph{Left}: $w=70$. \emph{Middle}: The eigenvalues with $w=100$ (red) added to the plot with the eigenvalues with $w=100$ (blue). \emph{Right}: Zoomed in the middle figure. The most outer blue and red eigenvalues agree to two or three digits. \label{FIG:vort-eigs-70-100}}
 \end{figure}

\begin{figure}[h]	
\includegraphics[width=0.34\linewidth,height=1.5in]{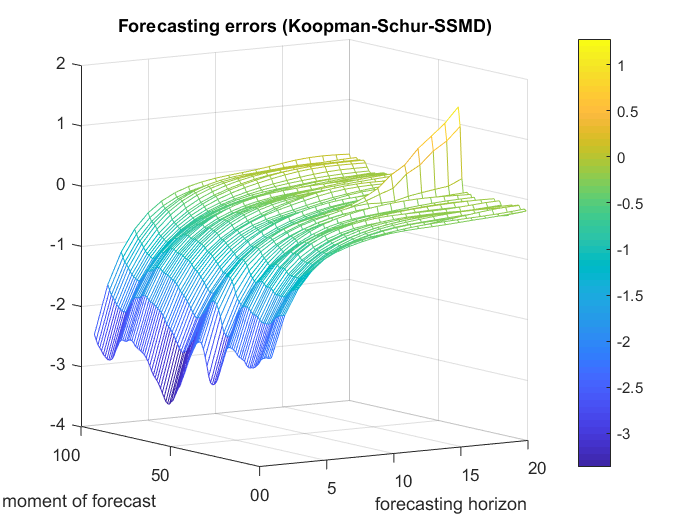}\hspace{-4mm}
\includegraphics[width=0.34\linewidth,height=1.5in]{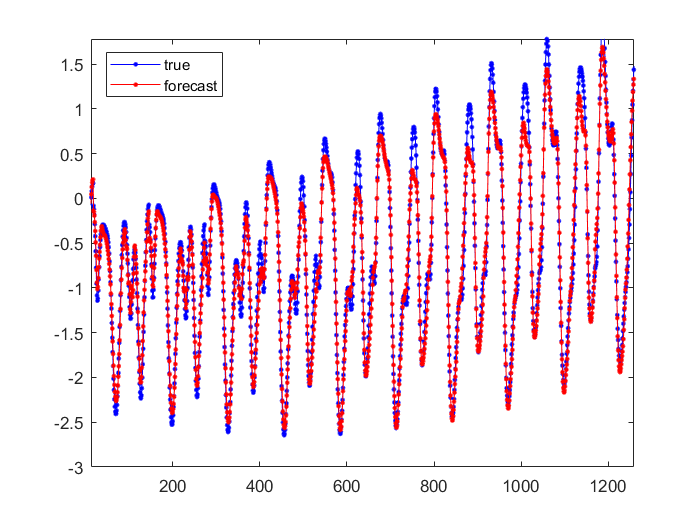}\hspace{-4mm}
\includegraphics[width=0.34\linewidth,height=1.5in]{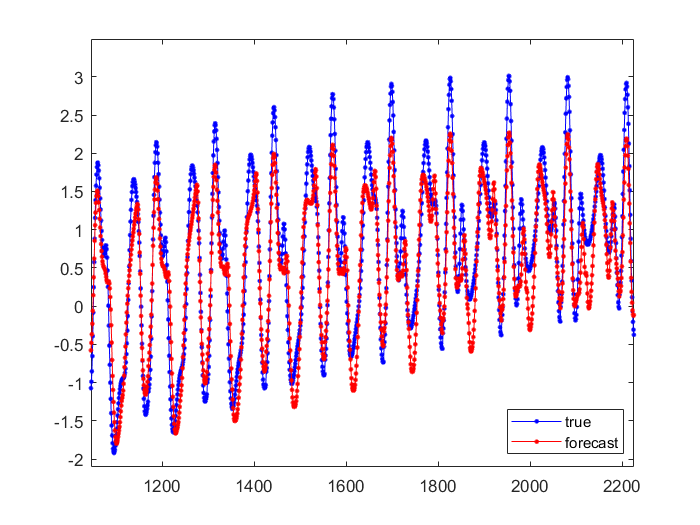}
\caption{(Example \ref{EX:vorticity}) \emph{Left}: The relative errors (in norm) with the prediction horizon $\tau=20$. \emph{Middle}: A detail of one forecast for $5$ steps ahead, zoomed to a part of the snapshot. 
\emph{Right}: A detail of one forecast for $10$ steps ahead, zoomed to a part of the snapshot.
}
 \end{figure}
 }
\end{example}

\section{Conclusions and ongoing and future work}\label{S=Conclusion}
We have presented a new theoretical and computational framework for Koopman operator based data driven analysis of nonlinear dynamics.  Preliminary testing has shown good performance and, most importantly, the KS-SSMD is based on the unitary Schur decomposition and it is thus independent of the diagonalizability of the EDMD matrix, or the (potentially high) condition number of the eigenvectors.

The ongoing and planned work, based on this report, includes 
\begin{itemize}
\item Further operator theoretic analysis of the approximation error.
\item A LAPACK style implementation of the presented algorithms, following our recent implementations of the DMD and the Hermitian DMD \cite{lawn-298-drmac}, \cite{drmac-DMD-2023}. (The Hermitian DMD \cite{lawn-300-drmac}, \cite{drmac-SYDMD-2023} is a special case of the KS-SSMD, where the Schur form is actually diagonal.)
\item A QR compressed implementation, including a streaming KS-SSMD based on fast updating/down-dating of the QR compressed representation of the snapshots.
\item Algorithmic details and software implementation for the applications:
\begin{itemize}
\item Identification of the coherent structures and their representations using flags of invariant subspaces.
\item Forecasting with adaptation to rapid changes in the dynamics.
\end{itemize}
\end{itemize}

\newpage

\bibliography{references}
\bibliographystyle{plain}
\end{document}